\documentclass[reqno,a4paper,12pt]{amsart} 

\usepackage{amsmath,amscd,amsfonts,amssymb}
\usepackage{mathrsfs,dsfont}

\usepackage{tikz}
\usetikzlibrary{patterns}
\usepackage{float}
\usepackage{appendix}
\usepackage{caption}
\usepackage{subcaption}
\usepackage{enumerate}
\usepackage{hyperref}

\numberwithin{equation}{section}
\numberwithin{figure}{section}

\def\R{\mathbb{R}}

\def\Z{\mathbb{Z}}

\def\1{\mathds{1}}

\def\dH{\dim_{\mathcal{H}}}

\renewcommand\leq{\leqslant}
\renewcommand\geq{\geqslant}

\renewcommand\hat{\widehat}

\newcommand{\supp}{\operatorname{supp}}

\newcommand{\dist}{\operatorname{dist}}

\theoremstyle{plain}
\newtheorem{thm}{Theorem}[section]
\newtheorem{lem}[thm]{Lemma}

\newtheorem{prop}[thm]{Proposition}

\newtheorem*{claim*}{Claim}
\newtheorem*{thm*}{Theorem}

\theoremstyle{definition}

\newtheorem*{definition*}{Definition}
\newtheorem*{remarks*}{Remarks}
\newtheorem*{remark*}{Remark}

\newenvironment{enumerate-math}
{\begin{enumerate}
\addtolength{\itemsep}{5pt}
}
{\end{enumerate}}

\newenvironment{enumerate-text}
{\begin{enumerate}
\addtolength{\itemsep}{5pt}
}
{\end{enumerate}}

\begin{document}

\title{Fourier Frames on Salem Measures}

\address{Department of Mathematics \& International Center for Mathematics, Southern
University of Science and Technology, Shenzhen, 518055, PR China}
\author{Longhui Li}
\email{12231267@mail.sustech.edu.cn}
\author{Bochen Liu}
\email{Bochen.Liu1989@gmail.com}

%\thanks{}
%\subjclass[2010]{}
\date{}

%\keywords{}

\begin{abstract}
For every $0<s\leq 1$ we construct $s$-dimensional Salem measures in the unit interval that do not admit any Fourier frame. Our examples are generic for each $s$, including all existing types of Salem measures in the literature: random Cantor sets (convolutions, non-convolutions), random images, and deterministic constructions on Diophantine approximations. They even appear almost surely as Brownian images. We also develop different approaches to prove the nonexistence of Fourier frames on different constructions. Both the criteria and ideas behind the constructions are expected to work in higher dimensions.

On the other hand, we observe that a weighted arc in the plane can be a $1$-dimensional Salem measure with orthonormal basis of exponentials. This leaves whether there exist Salem measures in the real line with Fourier frames or even orthonormal basis of exponentials a subtle problem.
\end{abstract}
\maketitle
%\setcounter{tocdepth}{1}
%\tableofcontents

\section{Introduction}
Suppose $\mu$ is a finite Borel measure on $\R^d$. A discrete set $\Lambda\subset\R^d$ is called a frame
spectrum of $\mu$ if there exist $0<A\leq B<\infty$ such that
$$A\|f\|_{L^2(\mu)}^2\leq \sum_{\lambda\in\Lambda}|\hat{f\,d\mu}(\lambda)|^2\leq B\|f\|_{L^2(\mu)}^2,\ \forall\,f\in L^2(\mu).$$
In this case we say $\mu$ admits a Fourier frame, or $\mu$ is frame-spectral, and call $A, B$ frame constants. When $A=B=1$ this inequality becomes the Parseval identity, and $\{e^{-2\pi i x\cdot\lambda}\}_{\lambda\in\Lambda}$ becomes an orthonormal basis of exponentials in $L^2(\mu)$. In this case we say $\mu$ is spectral and call $\Lambda$ a spectrum of $\mu$. In particular every spectrum is a frame spectrum.

Fourier frame was introduced by Duffin and Schaeffer \cite{DS52} in 1952, and has become a powerful tool in signal transmission and reconstruction, image processing and compression, data recovery and inpainting, etc. Though the Parseval identity is weakened to an inequality, some nice properties of orthonormal basis still remain on Fourier frames. For example every $f\in L^2(\mu)$ still has a Fourier series expansion. Although such an expansion may not be unique, the redundancy of frames provides extra advantages compared with orthonormal basis. We refer to \cite{You01}\cite{Christensen16} for basic properties of frames.

Fourier frames on sets (i.e. $d\mu=\chi_Edx$) is usually better understood than general measures. In \cite{NOU16}, Nitzan, Olevskii, and Ulanovskii proved that any Euclidean subset of finite Lebesgue measure (not necessarily bounded) admits Fourier frames. Because of this strong result, the discussion on the existence of Fourier frames is mainly on singular measures. In \cite{LW17}, Lai and Wang  constructed a family of singular measures that admit Fourier frames but no orthonormal basis of exponentials.

It is proved by Jorgensen and Pedersen \cite{JP98} that the standard one-fourth Cantor measure admits orthonormal basis of exponentials, while the standard one-third Cantor measure does not. It is a famous open problem in this area (proposed by Strichartz in  \cite{Strichartz00}) that whether the standard one-third Cantor measure admits Fourier frames.

In addition to finding Fourier frames, one can also study necessary conditions for its existence, or figure out which measures do not admit Fourier frames. In \cite{HLL13}, He, Lau and Lai proved that a frame-spectral measure must be of pure type, that is, it must be purely discrete, purely singular continuous, or purely absolutely continuous. So one can construct nonexistence examples by combining measures of different types. In \cite{DL14}, Dutkay and Lai obtained a family of nonexistence examples by analyzing a uniformity condition that generalizes \cite{Lai11} (see Section \ref{subsec-criteria} below). In particular self-affine measures under the non-overlapping condition with unequal weights do not admit any Fourier frame. See also \cite{FL18} for related examples.

In \cite{Lev18}, Lev studied several existence and nonexistence examples. We shall discuss his criterion for the nonexistence in Section \ref{subsec-criteria} below. In particular, he pointed out that a small spherical cap admits Fourier frames, and then asked about the surface measure on the whole sphere. Lev's question was answered negatively by Iosevich, Lai, Wyman, and the second author in \cite{ILLW19}. More generally, the surface measure on the boundary of any convex body with nonvanishing Gaussian curvature does not admit Fourier frames. Compared with other nonexistence examples, this is the only Salem measure (see Section \ref{subsec-Salem} below for a brief review on Salem sets and measures). Due to the lack of curvature in the real line, it cannot be extended to $\R$ immediately. In this paper we find its analogs in the real line. In fact there are many of them.

\begin{thm}
	\label{thm-main}
	For every $0<s\leq 1$ there exist $s$-dimensional Salem measures on the unit interval that do not admit any Fourier frame.
\end{thm}
If one only needs a single example, then a one-line construction based on any given Salem measure would work. See \eqref{one-line} below. However, the purpose of this paper is to show that Salem measures without Fourier frames are generic: they come out naturally from all existing constructions of Salem measures of different types, including random Cantor sets (convolutions (Theorem \ref{thm-Salem}) , non-convolutions (Theorem \ref{thm-non-convolution})), random images (Theorem \ref{thm-Brownian}), and deterministic constructions on Diophantine approximations (Theorem \ref{thm-Diophantine}). They can even appear almost surely as Brownian images (Theorem \ref{thm-Brownian}). We also develop different methods on different constructions to prove the nonexistence of Fourier frames. See Section \ref{subsec-criteria} below for a brief summary of our criteria. In this paper we only focus on examples in the real line, but both the criteria and ideas behind the constructions are expected to work in higher dimensions. As a remark, our result shows that, despite Fourier frame is a relaxation of orthonormal basis and exist on every set of finite Lebesgue measure, its existence on singular measures is not generally guaranteed.

After seeing so many nonexistence examples, one may wonder if no Salem measure admits Fourier frame. The answer is, no, at least in the plane. If we consider a measure on a small arc in the unit circle, say
$$\int f\,d\sigma:=\int^{\frac{1}{2}}_{-\frac{1}{2}}f(x, \sqrt{1-x^2})\,dx,$$
then clearly it has orthonormal basis of exponentials with $\Z\times\{0\}$ as a spectrum. Although its Fourier transform is not a standard stationary phase estimate (in which the integrand is $C_0^\infty$), the stationary phase technique still works and gives
$$|\widehat{\sigma}(\xi)|\lesssim |\xi|^{-\frac{1}{2}}.$$
See Section \ref{sec-arc} for computation. Therefore $\sigma$ is a $1$-dimensional Salem measure in the plane with orthonormal basis of exponentials. This makes whether there exist Salem measures in the real line with Fourier frames or even orthonormal basis of exponentials a subtle problem.

We would like to also point out that the cap example is more complicated in higher dimensions. For example the measure defined by
$$\int f\,d\sigma:=\int_{[-\frac{1}{2}, \frac{1}{2}]^{d-1}}f(x_1,\dots,x_{d-1}, \sqrt{1-|x|^2})\,dx,$$
though still spectral, may not be Salem. This is because 
$$\widehat{\sigma}(\xi_1,\dots,\xi_{d-1},0)=\prod_{i=1}^{d-1}\hat{\chi_{[-\frac{1}{2}, \frac{1}{2}]}}(\xi_i) $$
that does not have desired decay on $|\xi|$ when $d$ is large. It may become Salem if the unit square is replaced by other domains, but it is not the main focus of this paper. The planar arc is already enough for our purpose here.

\medskip

\subsection*{Organization}This paper is organized as follows. In Section \ref{sec-prelim} we give a brief review on the history of Salem sets and measures, then introduce criteria we shall use on the nonexistence of Fourier frames. The one-line construction of Salem measures without Fourier frames is also given in this section. The rest of this paper is on the constructions of Salem measures of different types and the proof on the nonexistence of Fourier frames : random Cantor sets as convolutions (Section \ref{sec-Cantor-convolution}), random Cantor sets as non-convolutions (Section \ref{sec-Cantor-non-convolution}), Brownian images (Section \ref{sec-random-images}), Diophantine approximations (Section \ref{sec-mu}-\ref{sec-final-proof}). This order is arranged based on the complexity of the construction and the proof:
$$\text{one-line}<\text{convolution}<\text{non-convolution}<\text{Brownian}<\text{Diophantine}.$$
In Section \ref{sec-arc} we discuss the Fourier decay of the weighted arc length measure on the unit circle.

\subsection*{Notation}
\ 

Fourier transform is defined by $\hat{f}(\xi)=\int_{\R} e^{-2\pi i x\cdot\xi}f(x)\,dx$ for $f\in L^1(\R)$ and $\hat{\mu}(\xi)=\int_{\R} e^{-2\pi i x\cdot\xi}d\mu(x)$ for finite Borel measures.
Fourier coefficient is defined by $\hat{F}(k)=\int_{-\frac{1}{2}}^{\frac{1}{2}} e^{-2\pi i x\cdot k}F(x)\,dx$ for $1$-periodic function $F$.

$X\lesssim Y$ means $X\leq CY$ for some constant $C>0$, and $X\lesssim_\epsilon Y$ means this constant $C$ may depend on $\epsilon$.

Throughout this paper $\phi_0\in C_0^\infty(-1,1)$ is a fixed even function satisfying $\phi_0, \hat{\phi_0}\geq 0$ and $\phi_0\geq \frac{1}{2}$ on $[-\frac{1}{2}, \frac{1}{2}]$. Such a function exists by taking $\phi_0=\varphi*\varphi$, where $\varphi\in C^\infty_0(-\frac{1}{2}, \frac{1}{2})$ is an arbitrary nonnegative even function satisfying $\varphi\geq 1$ on $[-\frac{1}{2}, \frac{1}{2}]$.

\subsection*{Acknowledgments} This work was partially supported by the National Key R\&D Program of China 2024YFA1015400, and the National Natural Science Foundation of China grant 12131011. We would like to thank Chun-Kit Lai for suggestions on the manuscript.

\section{Preliminaries}\label{sec-prelim}
\subsection{Salem sets and Salem measures}\label{subsec-Salem}
Hausdorff dimension of sets (denoted by $\dH$ in this paper) is originally defined by coverings. There are also equivalent definitions in terms of measures. Let $\mathcal{M}(E)$ denote the collection of nonzero finite Borel measures supported on a compact subset of $E\subset\R^d$. It is well known (see, e.g. Section 2.5 in \cite{Mat15}) that
\begin{equation}\label{def-Hausdorff-dim}\dH E=%&\sup\{s: \exists\,\mu\in\mathcal{M}(E), s.t.\,\mu(B(x,r))\leq r^s,\,\forall\,x\in\R^d,\,\forall\,r>0\}\\=&
\sup\{s: \exists\,\mu\in\mathcal{M}(E), s.t.\,I_s(\mu)<\infty\},\end{equation}
where
$$I_s(\mu):=\iint |x-y|^{-s}d\mu(x)d\mu(y)=C_{d,s}\int_{\R^d}|\hat{\mu}(\xi)|^2|\xi|^{-d+s}\,d\xi$$
denotes the $s$-dimensional energy. By this equivalent definition, if one defines the Fourier dimension by
    $$\dim_{\mathcal{F}} E:=\sup\{t\leq d:\exists\, \mu\in\mathcal{M}(E), s.t.\, |\hat{\mu}(\xi)|\lesssim |\xi|^{-t/2}\},$$
then $$\dim_{\mathcal{F}} E\leq \dH E.$$ A set is called Salem if the equality holds. A measure is called Salem if $$|\hat{\mu}(\xi)|\lesssim_\epsilon |\xi|^{-\dH(\supp\mu)/2+\epsilon}.$$

If $\mu$ is a Salem measure and $\mu=\mu_1+\mu_2$ is a nontrivial decomposition, then $I_s(\mu_i)\leq I_s(\mu)$,
which by \eqref{def-Hausdorff-dim} implies $\dH(\supp\mu_i)\geq \dH(\supp\mu)$, and therefore
$$\dH(\supp\mu_1)=\dH(\supp\mu_2)=\dH(\supp\mu).$$
In particular every $s$-dimensional Salem measure ($0<s<1$) has pure type, thus the criterion of He-Lai-Lau \cite{HLL13} does not apply.

We would like to classify all existing Salem sets/measures in $\R^d$ into the following three groups:

\begin{itemize}
	\item Random Cantor sets. It dates to Salem's original construction \cite{Salem51} in 1951. See also \cite{Bluhm96}\cite{Bluhm99}\cite{LP09}\cite{HL13}\cite{Che16}\cite{HL16}\cite{Shmerkin17}\cite{SS18}.
	\item Random images. It dates back to Kahane's work on images of Brownian motions \cite{Kahane66_Brownian} in 1966 (we learn this from \cite{Mat15}, Chapter 12). See also \cite{Kahane66_Gaussian}\cite{Kahane85}\cite{SX06}\cite{KWX06}\cite{WX07}\cite{Ekstrom16}.
	\item Diophantine approximations. It dates back to Kaufman's construction \cite{Kau81} in 1981.  See also \cite{Bluhm98}\cite{Hambrook17}\cite{FH23}\cite{FHR25}\cite{LL24+}\cite{FHR25+}. Nowadays Kaufman-type construction is still the only way to construct deterministic Salem sets.
\end{itemize}
There is also a Baire category argument, due to K\"orner (see, e.g. \cite{Kor11}\footnote{There is a gap in this paper, as explained in a footnote in \cite{LL24+}, but it is fine for Salem measures.}), that helps find Salem measures satisfying extra properties (see, e.g. \cite{CS17}).

\subsection{Criteria on the nonexistence of Fourier frames}\label{subsec-criteria}
There are not many general criteria to determine whether a singular measure admits Fourier frames, and we just mentioned that the pure type criterion does not work on Salem measures. A widely used one is the uniformity criterion due to Dutkay and Lai \cite{DL14}. It states that $\mu$ admits no Fourier frame, if there exists a nontrivial restriction $\mu|_F$ that is absolutely continuous with respect to a translated copy of $\mu$, with the Radon-Nikodym derivative unbounded. With this criterion one can quickly construct a Salem measure without Fourier frames: let $\nu$ be an arbitrary Salem measure on $[0,1]$ and $\psi\in C_0^\infty$, $\psi>0$ $\nu$-a.e. with essential infimum $0$ on $\supp\nu$. Then it is easy to verify that
\begin{equation}\label{one-line}d\mu(x):=d\nu(x+1)+\psi(x)d\nu(x) \text{ on }[-1,1]\end{equation}
is still Salem, and $d\nu(x+1)$ is absolutely continuous with respect to $\psi(x)d\nu(x)$ up to translation with Radon-Nikodym derivative $\psi^{-1}|_{\supp\nu}$, unbounded.

As promised in the introduction, we shall explore more constructions, and on some of which Dutkay-Lai's uniformity criterion still helps. It works particularly well on measures constructed by infinite convolutions of discrete measures (e.g. self-affine measures) under the non-overlapping condition. What we need is the following version.
\begin{prop}\label{prop-uniformity:label}
Suppose $\{d_k\}_{k=1}^{\infty}$ is an integer sequence with $d_k\geq2$. For each $k=1,2,\dots$, suppose $0<x^{k}_1<x^k_2\dots<x^k_{d_k}<1$ are distinct real numbers, $l_k$ is a positive number satisfying 
\begin{align}
    l_k<\min\limits_{1\leq j\leq d_k-1}\{x^k_{j+1}-x^{k}_{j}\},\ l_k<1-x^k_{d_k},
\end{align}
$p_{j,k}$ are positive weights with $\sum_{j=1}^{d_k} p_{j,k}=1$, and
$$\nu_k:=\sum_{j=1}^{d_k} p_{j,k}\cdot\delta_{(\prod_{n=1}^{k-1}l_n)x^k_{j}},$$
with $\prod_{n=1}^0l_n:=1$ as convention. Then the probability measure
$$\mu:=\lim_{k\rightarrow\infty}\nu_1*\nu_2*\cdots*\nu_k$$
admits no Fourier frame if 
$$\lim_{k\rightarrow\infty}\frac{\max\limits_{j_1,\dots,j_k}p_{j_1,1}\cdots p_{j_k,k}}{\min\limits_{j_1,\dots,j_k}p_{j_1,1}\cdots p_{j_k,k}} =\infty.$$
\end{prop}
\begin{proof}[Proof of Proposition \ref{prop-uniformity:label}]
	 For readers' convenience we sketch the proof for this particular case. The convergence of $\mu$ is standard by the choice of $l_k$. For fixed $j_1,\dots, j_N$, $j_,',\dots,j_N'$, consider sets
	 $$F:=\sum_{k=1}^N (\prod_{n=1}^{k-1}l_n)x^k_{j_k}+\sum_{k>N}\supp\nu_k,\quad  F':=\sum_{k=1}^N (\prod_{n=1}^{k-1}l_n)x^k_{j'_k}+\sum_{k>N}\supp\nu_k.$$
	 Then the key observation is $(p_{j_1,1}\cdots p_{j_k,k})^{-1}\mu|_F$ and $(p_{j'_1,1}\cdots p_{j'_k,k})^{-1}\mu|_{F'}$ coincide up to a translation $\sum_{k=1}^N (\prod_{n=1}^{k-1}l_n)(x^k_{j_k}-x^k_{j'_k})$, which implies
	 $$(p_{j_1,1}\cdots p_{j_N,N})^{-1}|\widehat{\mu|_F}(\xi)|=(p_{j'_1,1}\cdots p_{j'_N,N})^{-1}|\widehat{\mu|_{F'}}(\xi)|,\ \forall\,\xi\in\R.$$
	 Therefore, if $\mu$ admits a frame-spectrum $\Lambda$ with frame constants $0<A\leq B<\infty$, then
	 $$\frac{A}{B}\cdot\left(\frac{p_{j_1,1}\cdots p_{j_N,N}}{p_{j'_1,1}\cdots p_{j'_N,N}}\right)^2 \leq \frac{\sum_{\lambda\in\Lambda}|\widehat{\mu|_F}(\lambda)|^2}{\sum_{\lambda\in\Lambda}|\widehat{\mu|_{F'}}(\lambda)|^2}=\frac{p_{j_1,1}\cdots p_{j_N,N}}{p_{j'_1,1}\cdots p_{j'_N,N}},$$
	 contradiction.
\end{proof}

We notice that Salem's original construction is an infinite convolution. In Section \ref{sec-Cantor-convolution} we shall modify his construction to make it satisfy Proposition \ref{prop-uniformity:label}, with Fourier decay preserved. Moreover this measure is almost Ahlfors-David regular, i.e.
$$r^{s+\epsilon}\lesssim_\epsilon \mu(B(x,r))\lesssim_\epsilon r^{s-\epsilon},\ \forall\, x\in\supp\mu,$$
thus different from the one-line construction \eqref{one-line} and all other constructions in this paper.

As far as we know, Salem's construction is the only infinite convolution among all existing Salem measures. So we need other criteria for other constructions. One criterion proposed in the literature is due to Lev (Theorem 3.3 in \cite{Lev18}), which states that $\mu$ on $\R^d$ admits no Fourier frame if there exist $\alpha, \beta\in(0, d]$, $\frac{1}{\alpha}-\frac{1}{\beta}>\frac{1}{d}$, such that $|\hat{f\,d\mu}(\xi)|\lesssim |\xi|^{-\beta/2}$ for some nonzero $f\in L^2(\mu)$ and
\begin{equation}\label{criteria-Lev}\liminf_{R\rightarrow\infty}\frac{1}{R^{d-\alpha}}\int_{|\xi|<R}|\hat{\mu}(\xi)|^2\,d\xi>0.\end{equation}
Another criterion is due to Iosevich, Lai, Wyman and the second author (Theorem 1.3, 1.4 in \cite{ILLW19}), which states that $\mu$ on $\R^d$ admits no Fourier frame if there exist $\gamma\in(0,d]$, $C>0$, such that $|\hat{\mu}(\xi)|\lesssim |\xi|^{-\gamma/2}$ and
\begin{equation}\label{criteria-ILLW}\inf_{\lambda\in\R^d}|\lambda|^{-\gamma}\int_{|\xi|\leq C}|\hat{\mu}(\lambda + \xi)|^2\,d\xi>0.\end{equation}
%The relation between $\alpha$ and $\beta$ can be improved to $\alpha<\beta$ by Theorem 1.3 in \cite{ILLW19} (together with Lemma 3.1 in \cite{Lev18}). 
Both criteria work well on surface measures in \cite{Lev18}\cite{ILLW19}. However, neither of \eqref{criteria-Lev}\eqref{criteria-ILLW} is easy to verify on general measures. In this paper we find the following criterion that helps us derive the nonexistence of Fourier frames on non-convolution random Cantor sets and random images.

\begin{prop}\label{prop-heavy-intervals:label}
Suppose $\mu$ is a finite Borel measure on $\R^d$ and there exist $0\leq \alpha<\beta\leq d$ such that
\begin{enumerate}[(i)]
	\item $|\hat{f\,d\mu}(\xi)|\lesssim |\xi|^{-\beta/2}$ for some nonzero $f\in L^2(\mu)$, and
	\item $\sup_x\mu(B(x,r))\gtrsim r^{\alpha}$, $\forall\,r>0$.
\end{enumerate}
Then $\mu$ admits no Fourier frame.
\end{prop}
\begin{proof}[Proof of Proposition \ref{prop-heavy-intervals:label}]
	Assume $\Lambda$ is a frame spectrum of $\mu$ with frame constants $0<A\leq B<\infty$. Although stated as $|\hat{\mu}(\xi)|$, the proof of Theorem 1.3 in \cite{ILLW19} still works on $|\hat{f\,d\mu}(\xi)|$, thus the Fourier decay assumption implies
	$$\sum_{\lambda\in\Lambda\backslash\{0\}}|\lambda|^{-\beta}=\infty.$$
	On the other hand, an argument of Shi in \cite{Shi21} says: since for each $x_0$ and each $|\xi|\leq (10r)^{-1}$,
	$$|\widehat{\mu|_{B(x_0,r)}}(\xi)|=\left|\int_{B(x_0,r)}e^{-2\pi i (x-x_0)\xi}d\mu\right|\geq \left|\int_{B(x_0,r)}\cos(2\pi(x-x_0)\xi) d\mu\right|\geq \frac{1}{2}\mu(B(x_0,r)),$$
	it follows that
	$$\#(\Lambda\cap B(0,(10r)^{-1}))\cdot\mu(B(x, r))^2\lesssim\sum_{\lambda\in \Lambda\cap B(0,(10r)^{-1})}|\widehat{\mu|_{B(x,r)}}(\lambda)|^2$$ $$\leq\sum_{\lambda\in \Lambda}|\widehat{\mu|_{B(x,r)}}(\lambda)|^2\\\leq B\cdot \mu(B(x,r)).$$
	Under our assumption this implies 
	$$\#(\Lambda\cap B(0,R))\lesssim R^\alpha,\ \forall\,R>0,$$
	contradiction to $\sum_{\lambda\in\Lambda\backslash\{0\}}|\lambda|^{-\beta}=\infty$ if $\alpha<\beta$. 
\end{proof}

A simple application of Proposition \ref{prop-heavy-intervals:label} is to consider $\mu=\mu_1+\mu_2$, where, for example, $\mu_1$ has Fourier decay of exponent $\beta/2$ and $\mu_2$ is $\alpha$-dimensional Ahlfors-David regular. But these sums can never be Salem measures as explained in the previous subsection. Then one may wonder if there exist Salem measures satisfying both conditions in Proposition \ref{prop-heavy-intervals:label} (with $f=1$), especially because people often study a set $E$ with Frostman measures satisfying
$$\mu(B(x,r))\lesssim r^{\dH E-\epsilon},$$
while $\dim_{\mathcal F}E\leq \dH E$. In fact for Salem measures it is possible to have ``heavy" intervals, and the only constrain for $\alpha, \beta$ is $\alpha\geq \beta/2$ due to
\begin{equation}\label{Mitsis:label}\begin{aligned}\mu(B(x,r))\lesssim \int \phi(\frac{y-x}{r})\,d\mu(y)=r\int e^{2\pi i x\cdot\xi}\,\overline{\hat{\phi}(r\xi)}\,\hat{\mu}(\xi)\,d\xi\\\leq r\int |\hat{\mu}(\xi)||\hat{\phi}(r\xi)|\,d\xi\lesssim r\int |\xi|^{-\beta/2}|\hat{\phi}(r\xi)|\,d\xi\lesssim r^{\beta/2},\end{aligned}\end{equation}
with $\phi\in C_0^\infty$ positive on $[-1,1]$. It is shown in \cite{LL24+} that for every $\alpha>\beta/2$ there are examples satisfying $\dH\supp\mu=\beta$,
$$|\hat{\mu}(\xi)|\lesssim |\xi|^{-\beta/2}, \quad \text{and}\quad  \mu(B(x_0,q_i^{-1}))\gtrsim q_i^{-\alpha}$$
for some $x_0\in\supp\mu$ and some rapidly increasing sequence $q_i$. But this is not enough for our use as $q_i$ has to increase rapidly due to technical reasons. In Section \ref{sec-Cantor-non-convolution}, \ref{sec-random-images} we shall see, given $\alpha>\beta/2$, there are many examples with $\dH\supp\mu=\beta$,
$$|\hat{\mu}(\xi)|\lesssim |\xi|^{-\beta/2}, \quad \text{and}\quad  \mu(B(x_0,r))\gtrsim r^{\alpha}, \ \forall\,r>0,$$
for some $x_0\in\supp\mu$. Such measures even appear almost surely as Brownian images.

\subsection{The key idea on Diophantine approximations}
The above covers all random constructions of Salem measures. It remains to deal with the only deterministic construction on Diophantine approximations, or more precisely, measures supported on
$$\bigcap_i\bigcup_{1\leq H\leq q_i^{s/2}}\mathcal{N}_{q_i^{-1}}\left(\frac{\Z}{H}\right).$$
For technical reasons, $q_i$ here has to be a rapidly increasing sequence. Therefore one can only have lower bound on a sequence of $q_i$-balls, not $r$-balls for arbitrary $r>0$. It is not an infinite convolution either. So no criterion above works. In fact we are not able to find any simple criterion that can help quickly determine the nonexistence of Fourier frames on Diophantine approximations. This makes this case much more complicated than the  others, and the proof takes nearly half of the pages in this paper (Section \ref{sec-mu}-\ref{sec-final-proof}). To give readers a hint, we are inspired by the argument on the sphere in \cite{ILLW19}, in which the key property on the spherical measure $\sigma$ is
\begin{equation}\label{property-sphere}R^{-d}\int_{|\xi|<R}|\hat{\sigma}(\lambda+\xi)|^2\,d\xi\approx |\lambda|^{-(d-1)},\end{equation}
uniformly in $|\lambda|>100$ and $10<R<|\lambda|/2$. This property highly relies on the asymptotic formula of the Fourier transform of the spherical measure, or more generally the Herz formula \cite{Herz62} for the boundary of convex bodies with nonvanishing Gaussian curvature. There seems no way to construct such a nice measure in the real line. The key idea in this paper is to construct a Salem measure $\mu$ with positive Fourier coefficients associated with an auxiliary measure $\nu\ll\mu$ such that
$$|\hat{\nu}(k+l)|\lesssim \hat{\mu}(k)+ (1+|k|)^{-0.99}$$
uniformly for, roughly speaking, $|l|\lesssim |k|$, $l,k\in\Z$. See Lemma \ref{lem-absolute-continuity}, \ref{lem-averaging-frequency} for precise statements. We believe our idea also works in higher dimensions. In particular we expect a variant of the very recent construction of Fraser-Hambrook-Ryou \cite{FHR25+} to be Salem in $\R^d$ without Fourier frames.

\section{Random Cantor sets: convolutions}\label{sec-Cantor-convolution}
In this Section we show that a minor modification of Salem's original construction is enough for applying Proposition \ref{prop-uniformity:label}, and the Fourier decay is preserved.
\begin{thm}\label{thm-Salem}
	For every $0<s\leq 1$, there exist a sequence of discrete probability measures $\nu_k=\sum\limits_{j:finite}p_{j,k}\delta_{x_{j,k}}$, such that
	$$\mu:=\lim\limits_{k\to\infty}\nu_1*\nu_2*\cdots\nu_k$$
	is a probability measure on $[0,1]$ with $\dH\supp\mu=s$, $|\hat{\mu}(\xi)|\lesssim_\epsilon|\xi|^{-s+\epsilon}$, while
	$$\lim_{k\rightarrow\infty}\frac{\max\limits_{j_1,\dots,j_k}p_{j_1,1}\cdots p_{j_k,k}}{\min\limits_{j_1,\dots,j_k}p_{j_1,1}\cdots p_{j_k,k}} =\infty.$$
	Moreover, $\mu$ is almost Ahlfors-David regular, namely
$$r^{s+\epsilon}\lesssim_\epsilon \mu(B(x,r))\lesssim_\epsilon r^{s-\epsilon},\ \forall\, x\in\supp\mu,$$
\end{thm}

 The key lemma in \cite{Salem51} is the following, and we shall use it without any change.
 \begin{lem}[Section 3 in \cite{Salem51}]\label{Salem-lemma}
     Let $F(\xi)=\sum_{j=1}^dp_j e^{-2\pi i x_j\xi}$, $p_j\in\mathbb{R}$. Given $r>0$, suppose
    \begin{align}
         A_r:=\min\limits_{\substack{m_j\in\mathbb{Z},\;|m_j|\leq 2r,\;1\leq j\leq d\\
    (m_1,\dots,m_d)\neq0\\\sum_{j=1}^dm_j=0}}|{\sum_{j=1}^dm_jx_j}|>0.
     \end{align}
     Then there exists $T_0>0$ depending on $r,d, p_j,x_j$ such that
     \begin{align}
         \frac{1}{T}\int_{a}^{a+T}|F(\xi)|^{2r}d\xi<2r^r(\sum_{j=1}^d p_j^2)^r,\quad \forall\,T\geq T_0,\ \forall\,a\in \R.
     \end{align}
 \end{lem}
The proof of Lemma \ref{Salem-lemma} is direct computation, from which one can easily see that $$T_0:=\frac{(|p_1|+\cdots|p_d|)^{2r}}{\pi\cdot A_r\cdot r^r(\sum_{j=1}^dp_j^2)^r}$$ is sufficient.

Now we state the construction. Let $\{d_k\}_{k=1}^{\infty}$ be an increasing sequence of integers with $d_k\geq2$ (we shall take $d_k=k+1$). Let $L_k$ be a sequence of positive numbers satisfying
\begin{equation}\label{choice-L-k}\frac{\log d_k}{\log(1/L_k)}=s,\end{equation}
and choose $d_k$ distinct numbers $0<x^{k}_1<x^k_2\dots<x^k_{d_k}<1$ such that
$$0<x^k_1<d_k^{-1}-L_k\quad \text{and}\quad   L_k<x^k_{j}-x^{k}_{j-1}<d_k^{-1}, \ \forall\,j\geq 2.$$
In particular $x_{d_k}^k<1-L_k$. So the condition in Proposition \ref{prop-uniformity:label} is satisfied for every $l_k\leq L_k$.
Moreover, by a simple volume estimate argument (see e.g. Lemma 7.3 in \cite{LP09}), for an increasing sequence $r_k\rightarrow\infty$ (we shall take $r_k=\sqrt{\log k}$), $x_j^k$ can be chosen to satisfy the extra condition
\begin{equation}\label{lower bound for linear combination}
         A_{r_k}:=\min\limits_{\substack{m_j\in\mathbb{Z},\;|m_j|\leq 2r_k,\;1\leq j\leq d_k\\
    (m_1,\dots,m_{d_k})\neq0\\\sum_{j=1}^{d_k}m_j=0}}|{\sum_{j=1}^{d_k}m_{j}x^k_{j}}|\geq (Cr_k)^{-2d_k},
  \end{equation}
that helps when applying Lemma \ref{Salem-lemma}.

Then, for each $t_k\in[0,1]$, define
$$l_k:=(1-\frac{1}{(k+1)^2})L_k+\frac{t_kL_k}{(k+1)^2}.$$ 
Notice that
\begin{equation}\label{range-l-k:label}(1-\frac{1}{(k+1)^2})L_k\leq l_k\leq L_k,\end{equation}
so $\prod_{k=1}^nl_k\approx\prod_{k=1}^nL_k$ uniformly in $n$.

Now we start from the unit interval $[0,1]$ and take 
\begin{align*}
    E_1:=\bigcup_{j_1=1}^{d_1}[x^1_{j_1}, x^1_{j_1}+l_1]=\bigcup_{j_1=1}^{d_1}I^1_{j_1}.
\end{align*}
For each $I_{j_1}$, let $a_{j_1}^1$ denote the left endpoint of $I_{j_1}$, and take
\begin{align*}
    E_2=\bigcup_{j_2=1}^{d_2}\bigcup_{j_1=1}^{d_1}[a_{j_1}^1+x^2_{j_2}l_1, a_{j_1}^1+x^2_{j_2}l_1+l_1l_2].
\end{align*}
Notice that each interval in $E_1$ is treated in the same way, which guarantees the measure defined below is a convolution.

Continue in this way and get $E_k$ for each $k\geq 1$, and let 
$$E_t:=\bigcap_{k=1}^{\infty}E_k,$$
then
\begin{equation}\label{E_t}
    E_t=\{\sum_{k=1}^{\infty}(\prod_{j=1}^{k-1}l_j)a_k:a_k\in\{x^k_1,x^k_2,\dots,x^k_{d_k}\},\;k\geq1\},
\end{equation}
with $\prod_{n=1}^0l_n:=1$ as convention. We use the notation $E_t$ since $E_t$ depends on the random variable $t=(t_k)_{k=1}^\infty\in[0,1]^\infty$.

There is a class of natural measures on $E_t$. For each $k$, choose a set of positive weights $\{p_{j,k}\}_{j=1}^{d_k}$ with $\sum_{j=1}^{d_k}p_{j,k}=1$. Then $E_t$ admits a natural probability measure $\mu_t$ defined by 
\begin{align}\label{infinite convolution}
\mu_t=\lim\limits_{k\to\infty}\nu_1*\nu_2*\cdots\nu_k,\quad \nu_k:=\sum_{j=1}^{d_k}p_{j,k}\cdot\delta_{(\prod_{n=1}^{k-1}l_n)x^k_{j}}.
\end{align}
The convergence is standard. We omit details. It is supported on $E_t$ because of \eqref{E_t}.

Here comes the only modification of us. In Salem's original construction $p_{j,k}=d_k^{-1}$ for each $j$, while we take
$$p_{1,k}=(1-\frac{1}{k+1})\frac{1}{d_k},\ p_{2,k}=(1+\frac{1}{k+1})\frac{1}{d_k},\ p_{j,k}=\frac{1}{d_k},\ j\geq 2. $$
Then
$$\lim_{k\rightarrow\infty}\frac{\max\limits_{j_1,\dots,j_k}p_{j_1,1}\cdots p_{j_k,k}}{\min\limits_{j_1,\dots,j_k}p_{j_1,1}\cdots p_{j_k,k}} =\prod_{k=1}^\infty\frac{(1+\frac{1}{k+1})}{(1-\frac{1}{k+1})}=\infty.$$

It remains to check its Fourier decay. Due to the convolution structure,
\begin{align}\label{formula of fourier transform}
    \widehat{\mu_t}(\xi)=\prod_{k=1}^{\infty}\widehat\nu_k
    =\prod_{k=1}^{\infty}P_k((\prod_{n=1}^{k-1}l_n)\xi),
\end{align}
where 
$$P_k(\xi):=\sum_{j=1}^{d_k}p_{j,k}e^{-2\pi ix^k_{j}\cdot \xi}.$$
Then for any positive integer $N$,
 \begin{equation}\label{Fourier-transform-mu-t}
     \int |\widehat{\mu_t}(m)|^{2r_N}dt\leq \int_0^1 \cdots\int_0^1|\prod_{k=1}^{N}P_{k+1}((\prod_{n=1}^{k}l_n)m)|^{2r_N}dt_1\cdots dt_{N}.
 \end{equation}
 Recall that $l_n=(1-\frac{1}{(n+1)^2})L_n+\frac{t_nL_n}{(n+1)^2}$, $t_n\in [0,1]$, then for $2\leq k+1\leq N+1$,
 \begin{align*}
     \int_0^1 |P_{k+1}((\prod_{n=1}^{k}l_n)m)|^{2r_{N}}dt_k&=\int_0^1 |P_{k+1}((\prod_{n=1}^{k-1}l_n)((1-\frac{1}{(k+1)^2})L_k+\frac{t_kL_k}{(k+1)^2})m)|^{2r_N}dt_k\\
     &=\int_0^1 |P_{k+1}(a_kt_k+b_k)|^{2r_N}dt_k\\
     &=\frac{1}{|a_k|}\int_{b_k}^{a_k+b_k}|P_{k+1}(\xi)|^{2r_N}d\xi
 \end{align*}
 for some $a_k$, $b_k$ satisfying $$|a_k|\geq C_0\frac{(\prod_{n=1}^kL_n)|m|}{(k+1)^2},$$ where $C_0>0$ is an absolute constant. To apply Lemma \ref{Salem-lemma}, we need
 \begin{equation}\label{greater-than-C-0}C_0\frac{(\prod_{n=1}^{N}L_n)|m|}{(N+1)^2}\geq T_0.\end{equation}
 By \eqref{lower bound for linear combination} and our choice of $p_{j,k}$,
 \begin{equation}\label{upper-bound-C-0}T_0:=\frac{(p_{1,N}+\cdots +p_{d_N,N})^{2r_N}}{\pi\cdot A_{r_N}\cdot r_N^{r_N}(\sum_{j=1}^{d_N}(p_{j,N})^2)^{r_N}}\leq \frac{(Cr_{N})^{2d_{N}}d_{N}^{r_{N}}}{\pi \cdot r_{N}^{r_{N}}}\cdot \left(\frac{1}{1+\frac{2}{(N+1)^2d_N}}\right)^{r_N}\leq \frac{(Cr_{N})^{2d_{N}}d_{N}^{r_{N}}}{\pi \cdot r_{N}^{r_{N}}}.\end{equation}
 Here, compared with Salem's original proof, there is an extra term $\left(\frac{1}{1+\frac{2}{(N+1)^2d_N}}\right)^{r_N}$ in the computation but it does not change the upper bound eventually used in the proof. Then by \eqref{greater-than-C-0}\eqref{upper-bound-C-0} and our choice of $L_k$ in \eqref{choice-L-k}, to apply Lemma \ref{Salem-lemma}, a condition
 \begin{equation}\label{modified condition for lemma}
     \log|m|\geq \sum_{n=1}^{N}s^{-1}\log d_n+2\log(N+1)+r_{N}\log d_{N}+2d_{N}\log (Cr_{N})-\log C_0
 \end{equation}
is sufficient. If we take $d_k=k+1$ and $r_k=\sqrt{\log k}$, we have $r_k^{d_k}>d_k^{r_k}$. So the condition
 \begin{align}\label{second modified condition}
     \log|m|\geq \sum_{n=1}^{N}s^{-1}\log d_n+2\log(N+1)+3d_{N}\log (Cr_{N})-\log C_0:=f(N)
 \end{align}
is also sufficient, which will be useful later.

As the right hand side of \eqref{second modified condition} is increasing in $N$, for each $m$ one can take $N_m$ to be the largest integer such that \eqref{second modified condition} holds. Then \eqref{second modified condition} holds for every $k\leq N_m$. Also
 \begin{align*}
     f(N_m)\leq\log|m|\leq f(N_m+1),
 \end{align*}
which implies
\begin{equation}\label{alpha-m}\alpha_m:=(s\log|m|)^{-1}\sum_{n=1}^{N_m}\log d_n\rightarrow 1,\ \text{as }m\rightarrow\infty,\end{equation}
because $d_{k+1}=o(\sum_{n=1}^k\log d_n)$ and $r_k=o(d_k)$ given $d_k=k+1$, $r_k=\sqrt{\log k}$.

Now we invoke Lemma \ref{Salem-lemma} to obtain, for each $k\leq N_m$, 
 \begin{align*}
      \int_0^1 |P_{k+1}((\prod_{n=1}^{k}l_n)m)|^{2r_{N_m}}dt_k&=\frac{1}{|a_k|}\int_{a_k/b_k}^{a_k/b_k+1/b_k}|P_{k+1}(\xi)|^{2r_{N_m}}d\xi\\
      &\leq(\frac{1}{|a_k|}\int_{a_k/b_k}^{a_k/b_k+1/b_k}|P_{k+1}(\xi)|^{2r_k}d\xi)^{\frac{r_{N_m}}{r_k}}\\
      &\leq 2^{r_N+1}r_k^{r_{N_m}}/d_k^{r_N},
 \end{align*}
where the difference from the Salem's computation is just an extra factor $2$ due to $\sum_{j=1}^{d_k}(p_{j,k})^2\leq\frac{2}{d_k}$ in place of $\sum_{j=1}^{d_k}(p_{j,k})^2=\frac{1}{d_k}$. Plugging this estimate into \eqref{Fourier-transform-mu-t}, we have
\begin{align}\label{integral estimate}
     \int |\widehat{\mu_t}(m)|^{2r_{N_m}}dt\leq \frac{2^{N_m(r_{N_m}+1)}\prod_{k=1}^{N_m}r_k^{r_{N_m}}}{\prod_{k=1}^{N_m}d_k^{r_{N_m}}}\leq \frac{2^{N_m(r_{N_m}+1)}r_{N_m}^{N_mr_{N_m}}}{\prod_{k=1}^{N_m}d_k^{r_{N_m}}}.
 \end{align}
By \eqref{alpha-m} and the fact $\sum_{n=1}^{N_m}\log d_n\sim N_m\log N_m$ (recall $d_k=k+1$), there exists $m_0$ such that
$$\int |\widehat{\mu_t}(m)|^{2r_{N_m}}dt\leq C_{m_0}|m|^{-\alpha_mr_{N_m}s},\ \forall\,|m|\geq m_0.$$
Hence
\begin{align}
    \int \sum_{m\in\Z}|m|^{\alpha_mr_{N_m}s-2}|\widehat{\mu_t}(m)|^{2r_{N_m}}dt<\infty,
\end{align}
which implies (recall $\alpha_m\rightarrow 1$ in \eqref{alpha-m})
\begin{align}
    |\widehat{\mu_t}(m)|\lesssim_\epsilon |m|^{-s/2+\epsilon},\ \forall\,m\in\Z,
\end{align}
for almost all $t\in[0,1]^\infty$. To see these measures indeed satisfy
\begin{equation}\label{Salem-mu-t}
    |\widehat{\mu_t}(\xi)|\lesssim_\epsilon |\xi|^{-s/2+\epsilon},\ \forall\,\xi\in\R\backslash\{0\},
\end{equation}
notice that 
$$\supp\mu_t\subset[x_1^1, x_{d_1}^1+L_1]\subset (0,1).$$
Therefore $d\mu_t=\psi\,d\mu_t$ for some $\psi\in C_0^\infty(0,1)$ with $\psi=1$ on $\supp\mu_t$. Then \eqref{Salem-mu-t} holds by standard arguments (see, e.g. Lemma 9.A.4 in \cite{Wol03}).

It is easy to see $$\dH \supp\mu\leq s$$ by trivial coverings, and \eqref{Salem-mu-t} gives $$\dH\supp\mu_t\geq \dim_{\mathcal F}\supp\mu_t\geq s,$$ thus $\mu_t$ is a $s$-dimensional Salem measure.

In the end we show this measure is almost Ahlfors-David regular. If $B(x,r)$ is an interval of the set obtained in the $n$th step of the construction of $E$, then  it follows directly from our construction that
$$(n+1)^{-1}r^s\leq\mu(B(x,r))\leq (n+1)r^s.$$
We call these intervals $n$-intervals. For general $B(x,r)$ centered at $x\in\supp\mu$, choose $n$ such that 
$$\prod_{k=1}^nl_k\leq r\leq \prod_{k=1}^{n-1}l_k.$$
Then $B(x,r)$ contains at least one $n$-interval but intersects at most two $(n-1)$-intervals. Hence
\begin{align*}
    r^{s+\epsilon}\lesssim_\epsilon (n+1)^{-1}(\prod_{k=1}^nl_k)^s\leq \mu(B(x,r))\leq 2n(\prod_{k=1}^{n-1}l_k)^s\lesssim_\epsilon r^{s-\epsilon}
\end{align*}
by our choice of $l_k$ (recall \eqref{choice-L-k}\eqref{range-l-k:label} and $d_k=k+1$).

\section{Random Cantor sets: non-convolutions}\label{sec-Cantor-non-convolution}

As mentioned in Section \ref{sec-prelim}, most Salem measures on random Cantor sets are not infinite convolutions. In this section we take Chen's construction \cite{Che16} (inspired by Laba-Praminik \cite{LP09}) as the model case and modify it to be a Salem measure satisfying Proposition \ref{prop-heavy-intervals:label}, thus admit no Fourier frames. The same idea should work on other constructions in this class.
\begin{thm}\label{thm-non-convolution}
	Given $0<s\leq 1$, there exist a sequence of dyadic nodes $A_{j}\subset\frac{\{0,1,\dots,2^{j}-1\}}{2^{j}}$ and a probability measure $\mu$
whose support is the decreasing limit
	$$E:=\bigcap_{j=1}^\infty\bigcup_{a\in A_j}[a,a+2^{-j}],$$
 such that $\supp\mu=\dH E=s$,
	$$|\hat{\mu}(\xi)|\lesssim_\epsilon|\xi|^{-s/2+\epsilon}, \quad \text{and}\quad \mu(x_0,r)\gtrsim r^{s/2}, \,\forall r>0$$
	for some $x_0\in\supp\mu$.
\end{thm}

Now we start our construction. For $0<s\leq 1$, let $\{t_j\}$ be a $\{1,2\}$-sequence with
\begin{equation}\label{prod-t-j}
    2^{sj}\leq t_1\cdots t_j\leq 2^{sj+1}.
\end{equation}
For example one can take $t_1=2$,
$$t_j=\begin{cases}1, & \text{ if } t_1\dots t_{j-1} \geq 2^{sj}\\2, & \text{ otherwise }
\end{cases}.$$
It satisfies \eqref{prod-t-j} by induction.

Let $\{A_j\}_{j\geq0}$ be a sequence of sets with $A_{0}=\{0\}$, and once $$A_{j}\subset\frac{\{0,1,\dots,2^{j}-1\}}{2^{j}}$$ is determined, choose $t_{j+1}$ elements 
$$A_{j+1,a}\subset \frac{\{0,1\}}{2^{j+1}}$$
randomly for each $a\in A_j$ under evenly distributed probability. More precisely, when $t_{j+1}=1$, take $A_{j+1,a}$ to be $\frac{\{0\}}{2^{j+1}}$ or $\frac{\{1\}}{2^{j+1}}$ with probability half-half, independently for each $a\in A_j$; when $t_{j+1}=2$, take $A_{j+1,a}=\frac{\{0,1\}}{2^{j+1}}$ for all $a\in A_j$. Then denote 
$$A_{j+1}:=\bigcup_{a\in A_{j}}(a+A_{j+1,a})\subset\frac{\{0,1,\dots,2^{j+1}-1\}}{2^{j+1}}.$$

After $\{A_j\}$ is determined, let
$$E_j:=\bigcup_{a\in A_j}[a,a+2^{-j}]$$
and
\begin{equation}\label{phi-j-random-Cantor}\phi_j:=2^{j}\sum_{a\in A_j}p_{j,a}\mathbf1_{[a,a+2^{-j}]}.\end{equation}

We need to explain our choice of weights $\{p_{j,a}\}_{a\in A_j}$, which is different from Chen's evenly distributed weights. Take $p_{0,a}=1$ as convention.

If $t_{j+1}=1$, all weights are preserved from the previous step, namely 
\begin{align*}
    p_{j+1,a}=p_{j,a'},
\end{align*}
where $a'\in A_j$ is the unique element satisfying $a-a'\in\{0,1\}/2^{j+1}$.

If $t_{j+1}=2$, then $A_{j+1,a}=\{0,1\}/2^{j+1}$ for all $a\in A_j$. Let $a_j$ be the smallest element in $A_j\cap[\frac{1}{2},1]$ (well defined as $t_1=2$), and set
\begin{equation}\label{p-j+1-a-j}\begin{cases}
	p_{j+1,a}=\frac{1}{\sqrt{2}}p_{j,a},\;p_{j+1,a+2^{-(j+1)}}=(1-\frac{1}{\sqrt{2}})p_{j,a}, & \text { for }a=a_j\\
	p_{j+1,a}=p_{j+1,a+2^{-(j+1)}}=\frac{1}{2}p_{j,a}, & \text { for } a\neq a_j
\end{cases}.\end{equation}
In particular,
\begin{equation}\label{p-j-a-j}
    p_{j,a_j}=(t_1\cdots t_j)^{-1/2}\approx 2^{-sj/2},\;\forall j\geq1.
\end{equation}
As a remark, $\frac{1}{\sqrt{2}}$ is the largest possible number to ensure the resulting measure is $s$-dimensional Salem. One can see this from the argument below or compare between \eqref{p-j-a-j} and
\eqref{Mitsis:label}.

Finally take
\begin{align*}
    E=\bigcap_{j\geq1}E_j \text{ and } \mu=\lim\limits_{j\to\infty}\phi_j.
\end{align*}
The convergence is again standard. Notice that the randomness on each $a\in A_j$ is independent, so the resulting measure is not necessarily an infinite convolution.

By \eqref{p-j-a-j} and \eqref{prod-t-j}, for every $r>0$ there exists $x_r$ such that
$$\mu(B(x_r,r))\gtrsim r^{s/2}.$$
Moreover, by our choice of $a_j$, there exists $\delta>0$ ($\delta=2^{j_0}$ where $t_{j_0}$ is the second $t_{j_0}=2$) such that
$$\mu(B(x_r,r))\subset (\delta, 1-\delta),\ \forall\,r<\delta.$$
This means, to apply Proposition \ref{prop-heavy-intervals:label}, it suffices to show
\begin{equation}\label{random-Cantor-Fourier-decay}|\hat{\mu}(k)|\lesssim_\epsilon|k|^{-s/2+\epsilon},\ \forall\,k\in\Z.\end{equation}
Then similar to the previous section, $\psi\,d\mu$ is a $s$-dimensional Salem measure where $\psi\in C_0^\infty$, $\psi\geq 0$, $\psi=1$ on $[\delta, 1-\delta]$, and $\psi\,d\mu$ is an $s$-dimensional Salem measure without Fourier frames by Proposition \ref{prop-heavy-intervals:label}.

It remains to prove \eqref{random-Cantor-Fourier-decay}. We shall show that one can choose $A_j$ such that
\begin{align}\label{Fourier coefficients estimate 2}
    |\widehat{\phi_{j+1}}(k)-\widehat{\phi_j}(k)| \lesssim_\epsilon\min\{1,\frac{2^{j+1}}{|k|}\}\cdot2^{-(s-\epsilon)j/2},\;\forall k\in\Z,
\end{align}
with $\phi_0:=\chi_{[0,1]}$. This immediately implies \eqref{random-Cantor-Fourier-decay} by taking sum in $j$. The proof is an induction on $j$.

By \eqref{phi-j-random-Cantor},
\begin{align*}
    \widehat{\phi_j}(k)=\frac{1-e^{2\pi ik/2^j}}{2\pi ik/2^j}\sum_{a\in A_j}p_{j,a}e^{-2\pi ika},
\end{align*}
so $\widehat{\phi_{j+1}}(k)-\widehat{\phi_j}(k)$ equals
\begin{align*}
    \frac{1-e^{2\pi ik/2^{j+1}}}{2\pi ik/2^{j+1}}\left(\sum_{a\in A_{j}}e^{-2\pi i ka}\left(\sum_{b\in A_{j+1,a}}p_{j+1,a+b}e^{-2\pi ik b}-\frac{1}{2}\sum_{b\in\{0,1\}/2^{j+1}}p_{j,a}e^{-2\pi i k b}\right)\right).
\end{align*} 
The case $t_{j+1}=2$ is straightforward: by our choice of $p_{j,a}$ in \eqref{p-j+1-a-j}, only the term $a=a_j$ remains, so by \eqref{p-j-a-j},
\begin{align*}
    |\widehat{\phi_{j+1}}(k)-\widehat{\phi_j}(k)|\leq (\sqrt{2}-1)\left|\frac{1-e^{2\pi ik/2^{j+1}}}{2\pi ik/2^{j+1}}\right|\cdot p_{j,a_j}\lesssim\min\{1,\frac{2^{j+1}}{|k|}\}\cdot2^{-sj/2}.
\end{align*}

The case $t_{j+1}=1$ requires some probabilistic estimates as in \cite{Che16} (Lemma 1) and \cite{LP09} (Lemma 6.3).
\begin{lem}[Bernstein's inequality]\label{Bernstein's inequality}
    Suppose $X_1,\dots,X_n$ are independent random variables with $|X_i|\leq 1$, $\mathbb{E}X_i=0$ and $E|X_i|^2=\sigma_i^2$, $1\leq i\leq n$. Let $\sigma^2\geq\sum_{i=1}^n\sigma_i^2$, and assume that  $\sigma^2\geq 6n\lambda$. Then 
    \begin{align}
        \mathbb{P}\left(|\sum_{i=1}^nX_i|\geq n\lambda\right)\leq 4e^{-n^2\lambda^2/8\sigma^2}.
    \end{align}
\end{lem}

Recall $A_{j+1,a}$ consists of a single element $x(a)\in \{0,1\}/2^{j+1}$, randomly selected with probability $1/2$ on each value. For each $k$, let
\begin{align*}
    \chi_{a}(k)=p_{j,a}e^{-2\pi i ka}\left(e^{-2\pi ik x(a)}-\frac{1}{2}\sum_{b\in\{0,1\}/2^{j+1}}e^{-2\pi i k b}\right)
\end{align*}
be a random variable. It is easy to check that 
$$\mathbb{E}\chi_a(k)=0, \quad |\chi_a(k)|\leq 2p_{j,a},\quad \mathbb{E}|\chi_a(k)|^2\leq 4p_{j,a}^2.$$
Notice $\chi_a(k)$ is $2^{j+1}$-periodic, so it suffices to consider $k\in \{0,1,\dots,2^{j+1}-1\}$. 

Then we apply Lemma \ref{Bernstein's inequality} with 
 \begin{equation}\label{n-and-sigma}n=\#A_j=t_1\cdots t_j \text{ and } \sigma^2=\sum_{a\in A_j}p_{j,a}^2:=\sigma^2_j\end{equation} to obtain
\begin{align*}
    \mathbb{P}\left(\left|\sum_{a\in A_j}\chi_a(k)\right|\geq\lambda\right)\leq 4e^{-\lambda^2/8\sigma_j^2}.
\end{align*}
Taking sum in $k\in \{0,1,\dots,2^{j+1}-1\}$. It follows that
\begin{align*}
    \mathbb{P}\left(\left|\sum_{a\in A_j}\chi_a(k)\right|\leq\lambda,\;\forall k\in\Z\right)\geq 1-2^{j+3}e^{-\lambda^2/8\sigma_j^2},
\end{align*}
which is $>0$ if 
$$\lambda:=2\sqrt{2}\sigma_j\cdot\sqrt{\log(2^{j+4})}.$$
With this $\lambda$, it implies that one can choose $A_j$ such that, for all integer $k$,
\begin{align*}
   |\widehat{\phi_{j+1}}(k)-\widehat{\phi_j}(k)|=\left|\frac{1-e^{2\pi ik/2^{j+1}}}{2\pi ik/2^{j+1}}\right|\cdot\left|\sum_{a\in A_j}\chi_a(k)\right|\lesssim\min\{1,\frac{2^{j+1}}{|k|}\}\cdot \sigma_j\cdot\sqrt{\log(2^{j+4})}.
\end{align*}
Then to prove \eqref{Fourier coefficients estimate 2}, thus completes the proof, it remains to show
\begin{align}
    \sigma_j^2\lesssim_\epsilon 2^{-(s-\epsilon)j}.
\end{align}

Indeed, by the definition of $\sigma_j$ in \eqref{n-and-sigma} and $p_{j,a}$ in \eqref{p-j+1-a-j}\eqref{p-j-a-j}, we have the following recursive formula on $\sigma_j$ with $\sigma_0:=0$:
\begin{align*}
    \sigma_{j+1}^2=
    \begin{cases}
        \frac{1}{2}\sigma_j^2+(\frac{3}{2}-\sqrt{2})p_{j,a_j}^2\leq \frac{1}{2}\sigma_j^2+(t_1\cdots t_{j+1})^{-1},&\text{ if }t_{j+1}=2\\
        \sigma_j^2,&\text{ if }t_{j+1}=1
    \end{cases}.
\end{align*}
Then by induction
$$\sigma_j^2\leq \#\{k\leq j:t_k=2\}\cdot(t_1\cdots t_j)^{-1},$$
which is, by \eqref{p-j-a-j},
$$\lesssim_\epsilon 2^{-(s-\epsilon)j},$$
as desired.

\section{Random images}\label{sec-random-images}
In this section, we show images of the standard Brownian motion over some measure $\mu$ are Salem measures without Fourier frames almost surely. In this paper we only discuss measures in $\R$, but the argument also works in higher dimensions. The same idea should work on other random images as well.

Let $\omega(t)$ denote the standard Brownian motion, $\omega(0)=0$. It is well known that for every $0<\alpha<\frac{1}{2}$, almost all $\omega(t)$ is $\alpha$-H\"older continuous. We learn the proof from \cite{Fal14}, Chapter 16. So by considering coverings,
$$\dH\omega(E)\leq 2\dH E, \ a.s.$$
What Kahane proved is, given a measure $\mu$ on the unit interval with
\begin{equation}\label{Frostman-for-Brownian}\mu(B(x,r))\lesssim r^s,\end{equation}
then the push-forward measure $\omega_*\mu$, or denoted by $\mu_\omega$, satisfies
\begin{equation}\label{Fourier-decay-Brownian-image}|\widehat{\mu_\omega}(\xi)|\leq C(\omega, s) |\xi|^{-s}\log|\xi|,\ a.s.\end{equation}
We refer to \cite{Mat15}, Chapter 12, for a proof.

In this paper we show that the ball condition \eqref{Frostman-for-Brownian} is not necessary for the Fourier decay \eqref{Fourier-decay-Brownian-image}. The following is already sufficient, that has its own interest. 
\begin{thm}
	\label{thm-Brownian}
	Suppose $s>0$ and $\mu$ is a finite Borel measure on $[0,1]$ satisfying
	\begin{equation}\label{ball-near-0}
    \mu([0,r])\lesssim r^{s/2},\;\forall r>0,
    \end{equation}
    and
    \begin{equation}\label{ball-away-from-0}
    \mu([x,x+r])\lesssim \min\{1, x^{-s/2} r^s\},\;\forall x\in\supp\mu\backslash\{0\},\,\forall r>0.
\end{equation}
Then
$$|\widehat{\mu_\omega}(\xi)|\leq C(\omega, s) |\xi|^{-s}\log|\xi|,\ a.s.$$
\end{thm}
If, in addition, $\dH\supp\mu=s\in(0,\frac{1}{2}]$ and 
\begin{equation}\label{heavy-ball-near-0}\mu([0,r])\gtrsim r^{s/2},\end{equation}
then $\mu_\omega$ is a $2s$-dimensional Salem measure satisfying
$$\mu_\omega([0,r])\geq \mu([0,C_\omega r^{1/\alpha}))\gtrsim_{\omega} r^{s/(2\alpha)},\ a.s.$$
for every $0<\alpha<\frac{1}{2}$ by the H\"older continuity. Hence Proposition \ref{prop-heavy-intervals:label} applies and $\mu_\omega$ is a $2s$-dimensional Salem measure without Fourier frames almost surely.

\subsection{One example of $\mu$}
We first show that measures with $\dH\supp\mu=s$ satisfying \eqref{ball-near-0}\eqref{ball-away-from-0}\eqref{heavy-ball-near-0} exist. For readers' convenience we take several notation from the previous section to state the following example. 

Let $\{t_j\}$ be the sequence in  \eqref{prod-t-j}, and construct $\{A_j\}$ in a similar way without any randomness, that is, take $A_{j+1,a}=\{0\}$ whenever $t_{j+1}=1$. More precisely, now $A_0:=0$,
$$A_{j+1}:=\begin{cases}
	A_j+\frac{\{0,1\}}{2^{j+1}}, & \text{if }t_{j+1}=2\\A_j, & \text{if }t_{j+1}=1
\end{cases},$$
and
$$E:=\lim E_j, \text{ where } E_j:=\bigcup_{a\in A_j}[a, a+2^{-j}].$$
It is standard to show $\dH E=s$. We omit details.

Then we construct a measure $\mu$ on $E$ in a similar way as the previous section. In fact the only difference is to fix $a_j=0$. More precisely, take
\begin{equation}\label{phi-j-random-image}\phi_j:=2^{j}\sum_{a\in A_j}p_{j,a}\mathbf1_{[a,a+2^{-j}]} \quad \text{and}\quad  \mu:=\lim\phi_j,\end{equation}
where $p_{0,0}=1$, 
$$p_{j+1,a}=p_{j,a}, \text{ if } t_{j+1}=1,$$
 and
\begin{equation}\label{p-j+1-a-j-random-image}\begin{cases}
	p_{j+1,a}=\frac{1}{\sqrt{2}}p_{j,a},\;p_{j+1,a+2^{-(j+1)}}=(1-\frac{1}{\sqrt{2}})p_{j,a}, & \text {for } a= 0\\
	p_{j+1,a}=p_{j+1,a+2^{-(j+1)}}=\frac{1}{2}p_{j,a}, & \text {for } a\neq 0
\end{cases}, \text{ if }t_{j+1}=2.\end{equation}
In particular,
\begin{align*}
    p_{j,0}=(t_1\cdots t_j)^{-1/2}\approx 2^{-sj/2},
\end{align*}
which implies \eqref{ball-near-0}\eqref{heavy-ball-near-0}.

It remains to check \eqref{ball-away-from-0}. It suffices to consider dyadic numbers $r=2^{-j}$. 

When $0< x\leq 2r$, it follows directly from \eqref{ball-near-0} that
$$\mu([x,x+r])\leq\mu([0,3r])\lesssim r^{s/2}\lesssim x^{-s/2}r^s,$$
so from now we assume 
\begin{equation}\label{x>2r}x> 2r=2^{-(j-1)}.\end{equation}

As $s\leq \frac{1}{2}<1$, there must be infinitely many $t_j=1$, so  every $x\in E$ can be written uniquely as
$$x=\sum_{k\geq1}x_k2^{-k},\;x_k\in\{0,t_k-1\}.$$ 

Denote $k(x)$ as the smallest integer such that $x_{k(x)}=1$ (well defined as $x\neq 0$). Notice that $x\approx 2^{-k(x)}$, and $k(x)\leq j-1$ because $x> 2^{-(j-1)}$. Therefore, by our construction, 
\begin{align*}
    \mu([\sum_{k=1}^jx_k2^{-k},\sum_{k=1}^jx_k2^{-k}+2^{-j}])=
        (t_1\cdots t_{k(x)-1})^{-1/2}\cdot(1-\frac{1}{\sqrt{2}})\cdot (t_{k(x)+1}\cdots t_j)^{-1}
\end{align*}
\begin{equation}\label{measure-2-j-interval}=(\sqrt{2}-1)\cdot(t_1\cdots t_{k(x)})^{1/2}\cdot (t_1\cdots t_j)^{-1}\approx 2^{k(x)s/2}\cdot 2^{-sj}\approx x^{-s/2}\cdot 2^{-sj}.\end{equation}

Due to
\begin{equation}\label{2-j-approximation}0\leq x-\sum_{k=1}^jx_k2^{-k}\leq 2^{-j},\end{equation}
we notice that
\begin{equation}\label{observation-union}[x, x+2^{-j}]\subset \bigcup_{m=1}^2 [\sum_{k=1}^jx_k2^{-k}+(m-1)2^{-j},\sum_{k=1}^jx_k2^{-k}+m2^{-j}]\end{equation}
and each interval in the right hand side has nonzero measure only if its left endpoint lies in $E$. Also by \eqref{x>2r}\eqref{2-j-approximation},
$$\frac{x}{2}\geq 2^{-j}\geq x-\sum_{k=1}^jx_k2^{-k},$$
which implies
\begin{equation}\label{compare-left-end-point}x\lesssim  \sum_{k=1}^jx_k2^{-k}+(m-1)2^{-j},\ m=1,2.\end{equation}

Hence, by \eqref{measure-2-j-interval}\eqref{observation-union}\eqref{compare-left-end-point},
$$\mu([x,x+2^{-j}])\lesssim x^{-s/2}\cdot 2^{-sj},$$
as desired.

\subsection{Proof of Theorem \ref{thm-Brownian}}
Finally we prove Theorem \ref{thm-Brownian}. The reduction to the ball condition of $\mu$ is the same as the classical argument. See, e.g. Chapter 12 in \cite{Mat15}. We sketch it for completeness. Readers who are familiar with this step can jump to \eqref{expectation-hat-mu-omega-xi}. 

For each positive integer $q$,
\begin{align*}
    |\widehat{\mu_{\omega}}(\xi)|^{2q}=\int\cdots\int e^{2\pi i\xi\cdot(\sum_{j=1}^q\omega(t_j)-\sum_{j=1}^q\omega(u_j))}d\mu(t_1)\cdots d\mu(t_q)d\mu(u_1)\cdots d\mu(u_q)
\end{align*}
\begin{align*}
    =(q!)^2\int_{0<t_1<\cdots<t_q}\int_{0<u_1<\cdots<u_q}e^{2\pi i\xi\cdot(\sum_{j=1}^q\omega(t_j)-\sum_{j=1}^q\omega(u_j))}d\mu(t_1)\cdots d\mu(t_q)d\mu(u_1)\cdots d\mu(u_q).
\end{align*}
Write
\begin{align*}
    \sum_{j=1}^q\omega(t_j)-\sum_{j=1}^q\omega(u_j)=\sum_{j=1}^{2q}\epsilon_j\omega(v_j),
\end{align*}
where 
$$\epsilon_j=\begin{cases}
	1, & \text{if }v_j\in\{t_1,\dots,t_q\}\\
	-1, & \text{if } v_j\in\{u_1,\dots,u_q\}
\end{cases} $$
and $v_j$ is increasing. Since $\{\epsilon_j\}_{j=1}^{2q}\subset\{-1,1\}^{2q}$ with $\sum_{j}\epsilon_j=0$ uniquely determines the order of $t_1,\dots,t_q,u_1,\dots,u_q$, the integral above can be written as
\begin{align*}
    |\widehat{\mu_{\omega}}(\xi)|^{2q}=&(q!)^2\sum_{\substack{\{\epsilon_j\}_{j=1}^{2q}\subset\{-1,1\}^{2q}\\ \sum_{j}\epsilon_j=0}}\int_{0<v_1<\cdots<v_{2q}}e^{2\pi i\xi\cdot(\sum_{j=1}^{2q}\epsilon_j\omega(v_j))}d\mu(v_1)\cdots d\mu(v_{2q})\\=&(q!)^2\sum_{\substack{\{\epsilon_j\}_{j=1}^{2q}\subset\{-1,1\}^{2q}\\ \sum_{j}\epsilon_j=0}}\int_{0<v_1<\cdots<v_{2q}}e^{2\pi i\xi\cdot\sum_{j=1}^{2q}(\epsilon_j+\cdots+\epsilon_{2q})(\omega(v_j)-\omega(v_{j-1}))}d\mu(v_1)\cdots d\mu(v_{2q}),
\end{align*}
where $v_0:=0$ as convention (recall $\omega(0)=0$). 

Then we take the expectation of both sides. As
$$\mathbb{E}(e^{2\pi i\xi\cdot((\epsilon_j+\cdots+\epsilon_{2q})(\omega(v_j)-\omega(v_{j-1})))})= e^{-2\pi ^2|\xi|^2(\epsilon_j+\cdots+\epsilon_{2q})^2(v_j-v_{j-1})},$$
with $R_j:=2\pi^2 |\xi|^2(\epsilon_j+\cdots+\epsilon_{2q})^2$ we obtain
\begin{align*}
    \mathbb{E}(|\widehat{\mu_{\omega}}(\xi)|^{2q})=&(q!)^2\sum_{\substack{\{\epsilon_j\}_{j=1}^{2q}\subset\{-1,1\}^{2q}\\ \sum_{j}\epsilon_j=0}}\int_{0<v_1<\cdots<v_{2q}}e^{-\sum_{j=1}^{2q}R_j(v_j-v_{j-1})}d\mu(v_1)\cdots d\mu(v_{2q})\\\leq &(q!)^2\cdot\sum_{\substack{\{\epsilon_j\}_{j=1}^{2q}\subset\{-1,1\}^{2q}\\ \sum_{j}\epsilon_j=0}} \prod_{j\text{ even}}\iint_{v_j>v_{j-1}}e^{-R_j(v_j-v_{j-1})}\,d\mu(v_j)\,d\mu(v_{j-1}).
\end{align*}
The advantage of even $j$ is $R_j\geq|\xi|^2$ because $\epsilon_j+\cdots+\epsilon_{2q}\neq0$ (recall $\epsilon_j=\pm 1$). A simple combinatorial fact shows the number of choice of $\{\epsilon_j\}$ in the sum is $\frac{(2q)!}{(q!)^2}$. So we end up with 
\begin{equation}
	\label{expectation-hat-mu-omega-xi}
	\begin{aligned}\mathbb{E}(|\widehat{\mu_{\omega}}(\xi)|^{2q})\leq & (q!)^2\cdot \frac{(2q)!}{(q!)^2}\cdot\prod_{j\text{ even}}\iint_{v_j>v_{j-1}}e^{-|\xi|^2(v_j-v_{j-1})}\,d\mu(v_j)\,d\mu(v_{j-1})\\ = &(2q)!\cdot \left(\iint_{y>x}e^{-|\xi|^2(y-x)}\,d\mu(y)\,d\mu(x)\right)^q\\=&(2q)!\cdot \left(\int\int_0^1\mu([x, x+|\xi|^{-2}\log r^{-1}])\,dr\,d\mu(x)\right)^q.
	\end{aligned}	
\end{equation}

Now we use our ball conditions \eqref{ball-near-0}\eqref{ball-away-from-0} to study the integral in \eqref{expectation-hat-mu-omega-xi}. First, by \eqref{ball-away-from-0},
\begin{align*}\int\int_0^1\mu([x, x+|\xi|^{-2}\log r^{-1}])\,dr\,d\mu(x)\lesssim &\int\int_{0}^1\min\{1, x^{-s/2}|\xi|^{-2s}(\log r^{-1})^s\}\,dr\,d\mu(x)\\\lesssim &\int\min\{1, x^{-s/2}|\xi|^{-2s}\}\,d\mu(x)=I+II.\end{align*}
Then, by \eqref{ball-near-0},
$$I:=\int_{x\leq |\xi|^{-4}}d\mu(x)\lesssim |\xi|^{-2s}$$
and
$$II:=|\xi|^{-2s}\int_{x\geq |\xi|^{-4}}x^{-s/2}d\mu(x)= |\xi|^{-2s}\int_0^{|\xi|^{2s}}\mu([0,r^{-2/s}])\,dr$$
$$\lesssim |\xi|^{-2s}\left(1+\int_1^{|\xi|^{2s}}r^{-1}\,dr\right)\lesssim |\xi|^{-2s}\log|\xi|.$$
Overall,
\begin{align}\label{estimate of expectation}
    \mathbb{E}(|\widehat{\mu_{\omega}}(\xi)|^{2q})\lesssim (2q)!\cdot(|\xi|^{-2s}\log|\xi|)^q.
\end{align}

The rest is again the same as the classical argument. Choose a discrete set $Z\subset\R$ such that
$$\sum_{z\in Z}|z|^{-(s+2)}<\infty$$
and for every $|\xi|\geq 1$ there exists $z_\xi\in Z$ such that
$$|\xi-z_\xi|\leq |\xi|^{-s}.$$
One can take, for example,
\begin{align*}
 Z:=\bigcup_{k\geq0}\left(\{\xi\in\R:2^k\leq|\xi|<2^{k+1}\}\cap 2^{-(k+1)s}\Z\right).
\end{align*}
Then, with $q_z:=[\log |z|+1]$,
\begin{align*}
    \mathbb{E}\left(\sum_{z\in Z}|z|^{-(s+2)}\frac{|\widehat{\mu_{\omega}}(z)|^{2q_z}}{(2q_z)!\cdot(C(s)\log|z|\cdot|z|^{-2s})^{q_z}}\right)<\infty,
\end{align*}
which implies
\begin{align*}
    |\widehat{\mu_{\omega}}(z)|\leq C(\omega,s) |z|^{-s}\log|z|,\;\forall z\in Z,\ a.s.
\end{align*}
For general $|\xi|\geq1$, as $\widehat{\mu_{\omega}}$ is Lipschitz,
\begin{align*}
    |\widehat{\mu_{\omega}}(\xi)|\leq |\widehat{\mu_{\omega}}(z_\xi)|+|\widehat{\mu_{\omega}}(z_\xi)-\widehat{\mu_{\omega}}(\xi)|\leq|\widehat{\mu_{\omega}}(z_\xi)|+|\xi|^{-s}\leq C(\omega,s) |\xi|^{-s}\log|\xi|.
\end{align*}

\section{Diophantine approximations: the target measure $\mu$}\label{sec-mu}
\begin{thm}
	\label{thm-Diophantine}
	For each $0<s\leq 1$, there exist a rapidly increasing sequence $\{q_i\}$ and a finite Borel measure $\mu$ supported on 
	$$\bigcap_i\bigcup_{1\leq H\leq q_i^{s/2}}\mathcal{N}_{q_i^{-1}}\left(\frac{\Z}{H}\right)\cap[-\frac{1}{2}, \frac{1}{2}]$$
	with positive Fourier coefficients, $\dH\supp\mu=s$, $|\hat{\mu}(\xi)|\lesssim_\epsilon|\xi|^{-s/2+\epsilon}$, while $\mu$ does not admit any Fourier frame.
\end{thm}

Throughout this paper $\phi_0\in C_0^\infty(-1,1)$ is a fixed even function satisfying $\phi_0, \hat{\phi_0}\geq 0$ and $\phi_0\geq \frac{1}{2}$ on $[-\frac{1}{2}, \frac{1}{2}]$. Such a function exists by taking $\phi_0=\varphi*\varphi$, where $\varphi\in C^\infty_0(-\frac{1}{2}, \frac{1}{2})$ is an arbitrary nonnegative even function satisfying $\varphi\geq 1$ on $[-\frac{1}{2}, \frac{1}{2}]$.

\subsection{Kaufman's original construction}
Let $\|x\|:=\dist(x,\Z)$. It is a classical result from 1930s (\cite{Jar31}\cite{Bes34}) that the set on Diophantine approximation,
$$	\{x\in\R: \|qx\|\leq q^{1-\alpha}\ \text{for infinitely many integers }q\},\ \alpha>2,$$
has Hausdorff dimension $2/\alpha$. In 1981, Kaufman \cite{Kau81} showed that this set is a Salem set. One can also see Chapter 9 in \cite{Wol03} for a detailed proof. 

The original construction of Kaufman can be sketched as follows. Let $\phi\in C^2[-1,1]$ be an arbitrary nonnegative function, $M_i$ be a rapidly increasing integer sequence, $\Phi_i$ be the $1$-periodization of $M_i^{\alpha-1}\phi(M_i^{\alpha-1}\cdot)$ and $\Phi_{i,p}(x):=\Phi_i(px)$. Then take $\mu$ as a weak* limit point of the sequence of the  arithmetic mean of
$$\prod_{i=1}^nF_i:=\prod_{i=1}^n\frac{1}{\#\mathcal{P}_i}\sum_{p\in\mathcal{P}_i}\Phi_{i,p},$$
where $$\mathcal{P}_i:=\{M_i\leq p\leq 2M_i, prime\}.$$
We shall discuss two Kaufman-type constructions, and take advantages of their different properties.

\subsection{An auxiliary function $\phi$}
In Kaufman's original construction, and in fact all previous Kaufman-type constructions, an arbitrary nonnegative $\phi\in C^2(-1,1)$ is sufficient. But it is not enough for us. To construct our desired measure, we need an auxiliary function $\phi$ with the following properties:
\begin{enumerate}[(i)]
	\item $\phi\in C^2(\R)$;
	\item $\supp\phi\subset(-1,1)$;
	\item $\int\phi=1$;
	\item $\phi\geq 0$;
	\item $\hat{\phi}\geq 0$;
	\item $\hat{\phi}(\xi+r)\approx\hat{\phi}(\xi)$ uniformly in $\xi\in\R$ and $r\leq 1$.
\end{enumerate}

One can construct a such $\phi$ explicitly. For example it suffices to construct $\phi$ satisfying (ii)-(v) and
$$\hat{\phi}(\xi)\approx (1+|\xi|)^{-4},$$
that implies (i) and (vi). Now let 
$$\phi_1(x)=\chi_{[-1/2,1/2]}, \quad \phi_2(x)=2x|_{[-1/2,1/2]}.$$
Then
$$\hat{\phi_1}(\xi)=\frac{\sin\pi\xi}{\pi\xi},\quad  \hat{\phi_2}(\xi)=i\left(\frac{\cos\pi\xi}{2\pi\xi}-\frac{\sin\pi\xi}{2\pi^2\xi^2}\right)=i\frac{\pi\xi\cos\pi\xi-\sin\pi\xi}{\pi^2\xi^2}$$
are continuous functions. Notice that $\hat{\phi_1}(\xi), \hat{\phi_2}(\xi)$ cannot be both zero: $$|\hat{\phi_2}(\xi)|=|\pi\xi|^{-1}, \text{ if }\xi\in \Z\backslash\{0\}.$$ Then there are many ways to combine $\phi_0, \phi_1, \phi_2$ to obtain a desired $\phi$. For example, one can take
$$\phi(x):=A_1\phi_0(x)+A_2(\phi_1*\phi_1+\phi_2*\phi_2^-)*(\phi_1*\phi_1+\phi_2*\phi_2^-)(4x),$$
where $\phi_2^-(x):=\phi_2(-x)$ and $A_1, A_2>0$ are constants to ensure $\phi\geq 0$ and $\int\phi=1$.  Then
$$\hat{\phi}(\xi)=A_1\hat{\phi_0}(\xi)+\frac{A_2}{4}\left((\frac{\sin\pi\xi/4}{\pi\xi/4})^2+(\frac{\pi\xi/4\cos\pi\xi/4-\sin\pi\xi/4}{(\pi\xi/4)^2})^2\right)^2.$$
It is a continuous strictly positive function, and when $|\xi|\rightarrow\infty$,
$$\lim_{|\xi|\rightarrow\infty}(\pi \xi/4)^4\hat{\phi}(\xi)=\frac{1}{4}A_2>0,$$
which implies
$$\hat{\phi}(\xi)\approx (1+|\xi|)^{-4},$$
as desired.

In addition to explicit constructions in the physical space, one can also construct a such $\phi$ from the frequency side. According to the celebrated Paley-Wiener theorem (see, e.g. Section 2.4 in \cite{You01}), if $F(z)$ is an entire function of exponential type at most $1$ (i.e. $|F(z)|\lesssim e^{2\pi |z|}$) and $F|_{z\in\R}\in L^2(\R)$, then there exists $\phi\in L^2[-1, 1]$ such that
$$F(z)=\int\phi(x) e^{2\pi i zx}dx.$$
So we only need to find such an $F$ whose restriction to the real line satisfies (ii)-(v) and
$$F(\xi)\approx (1+|\xi|)^{-4}.$$
There are many options. For example one can take
$$F(z)=A_1\left(\frac{\pi z/4-\sin \pi z/4}{z^3}\right)^2+A_2\hat{\phi_0}(z),$$
where $A_1, A_2>0$ are properly chosen.

\subsection{Construction of the target measure $\mu$}
With this auxillary $\phi$, we can state our construction of $\mu$. Let $\{q_i\}_{i\geq 1}$ be a rapidly increasing sequence (not necessarily integer, $q_{i+1}\geq q_i^{10i}$ is sufficient). Denote
\begin{equation}\label{def-P-mu}\mathcal{P}^\mu_i:=\{1\}\cup\{p\leq q_i^{s/2}/2, prime\}.\end{equation}
For each $p\in \mathcal{P}^\mu_i$, let $\Phi_{i,p}^\mu$ be the $p^{-1}$-periodization of $p^{-1}q_i\phi(q_ix)$, namely 
\begin{equation}
\label{def-Phi-i-p-mu}
\Phi_{i,p}^\mu(x):=\sum_{v\in\mathbb{Z}}p^{-1}q_i\phi(q_i(x-\frac{v}{p}))=\sum_{n\in\Z}\hat{\phi}(pq_i^{-1}n) e^{2\pi i p n x},
\end{equation}
and take
\begin{equation}\label{def-F-i-mu}F^\mu_i:=\frac{1}{\#\mathcal{P}_i^\mu}\sum_{p\in\mathcal{P}^\mu_i}\Phi_{i,p}^\mu\end{equation}
which is $1$-periodic. Notice for every integer $k$,
\begin{equation}\label{Fourier-coefficient-F-i-mu}\widehat{F_i^\mu}(k)=\begin{cases}
	1, & k=0\\
	(\#\mathcal{P}_i^\mu)^{-1} \cdot \#\{p\in\mathcal{P}_i^\mu: k\in p\Z\}\cdot\hat{\phi}(q_i^{-1}k), &k\neq 0
\end{cases}.\end{equation}
In particular, for all integer $k\neq 0$,
\begin{equation}
	\label{F-i-mu-particular}
	\widehat{F_i^\mu}(k)\geq (\#\mathcal{P}_i^\mu)^{-1}\cdot\hat{\phi}(q_i^{-1}k)>0,
\end{equation}
that would imply the measure $\mu$ below has positive Fourier coefficients.

Finally we take $\mu$ as the weak* limit 
$$\lim_{n\rightarrow\infty}2\phi(2x)\prod_{i=1}^n F^\mu_i(x),$$
supported on
$$\bigcap_i\bigcup_{p\in\mathcal{P}_i^\mu}\mathcal{N}_{q_i^{-1}}\left(\frac{\Z}{p}\right)\cap[-\frac{1}{2}, \frac{1}{2}].$$
For convenience we denote $$F_0^\mu(x):=2\phi(2x).$$  This weak* convergence is guaranteed by an associated stability lemma like all Kaufman-type constructions.
\begin{lem}\label{stability-lemma}
	Suppose $\psi\in C^2([-\frac{1}{2}, \frac{1}{2}])$ and $\psi(-\frac{1}{2})=\psi(\frac{1}{2})$. Then
        $$|\widehat{\psi F_i^\mu}(k)-\widehat{\psi}(k)|\leq C (\|\psi\|_{L^\infty}+\|\psi''\|_{L^\infty})\cdot
        \begin{cases}
            q_i^{-s/2}(\log q_i)^2, & |k|\leq q_i\\
            |k|^{-s/2}(\log|k|)^2, & |k|\geq q_i
        \end{cases}.$$
\end{lem}
We leave the proof of Lemma \ref{stability-lemma} to the end of this section. As translated $\phi$ in the definition of $\Phi_{i,p}^\mu$ \eqref{def-Phi-i-p-mu} have disjoint support, 
	\begin{equation}\label{F-i-mu-L-infty}
	\|F_i^\mu\|_{L^\infty}\leq \sup_{p\in\mathcal{P}_i^\mu}\|\Phi_{i,p}^\mu\|_{L^\infty}\leq q_i\|\phi\|_{L^\infty}\end{equation}
	and
	\begin{equation}\label{F-i-mu''-L-infty} \|(F_i^\mu)''\|_{L^\infty}\leq \sup_{p\in\mathcal{P}_i^\mu}\|(\Phi_{i,p}^\mu)''\|_{L^\infty}\leq q_i^3\|\phi''\|_{L^\infty}.\end{equation}
Therefore
$$\|\prod_{i=1}^{n-1} F_i^\mu\|_{L^\infty}+\|(\prod_{i=1}^{n-1} F_i^\mu)''\|_{L^\infty}\leq n^2 C_\phi^{n-1}q_{n-1}^3\prod_{i=1}^{n-2} q_i.$$
Given $\{q_i\}$ a rapidly increasing sequence and $q_1$ large enough, by taking $\psi=\prod_{i=0}^{n-1} F^\mu_i$ it follows from Lemma \ref{stability-lemma} that
\begin{equation}\label{k=0}
	\sum_{n=1}^\infty |\widehat{\prod_{0\leq i\leq n}F_i^\mu}(0)-\widehat{\prod_{0\leq i\leq n-1}F_i^\mu}(0)|\leq C_{\epsilon}\sum_{i=1}^\infty q_i^{-s/2+\epsilon}<1
\end{equation}
and for all integer $k\neq 0$,
\begin{equation}\label{k-neq-0}\sum_{n=1}^\infty |\widehat{\prod_{0\leq i\leq n}F_i^\mu}(k)-\widehat{\prod_{0\leq i\leq n-1}F_i^\mu}(k)|\lesssim_\epsilon |k|^{-s/2+\epsilon}\cdot\#\{i:q_i\leq |k|\}+\sum_{i:q_i\geq |k|}q_i^{-s/2+\epsilon}\lesssim_\epsilon |k|^{-s/2+\epsilon}.	
\end{equation}
By \eqref{k=0}, there exists a subsequence of $\prod_{i=0}^n F_i^\mu$ that is weakly convergent to a nontrivial finite Borel measure $\mu$. Then the finiteness of \eqref{k-neq-0} implies that
$$\int \varphi \prod_{i=0}^n F_i^\mu=\sum_{k\in\Z}\widehat{\varphi}(k)\widehat{\prod_{0\leq i\leq n} F_i^\mu}(k)\rightarrow \int \varphi\,d\mu$$
for every smooth $1$-periodic function $\varphi$. Together one can conclude that
$$\prod_{i=0}^n F_i^\mu\xrightarrow{w^*}\mu.$$

Once the weak* convergence is verified, \eqref{k-neq-0} implies 
$$\left|\left(\chi_{[-\frac{1}{2}, \frac{1}{2}]}\prod_{i=1}^n F_i^\mu\right)^{\wedge}(k)\right|\lesssim_\epsilon|k|^{-s/2+\epsilon},$$
and therefore, after multiplying $2\phi(2x)$,
$$|\hat{\mu}(\xi)|\lesssim_\epsilon|\xi|^{-s/2+\epsilon}$$
follows by standard arguments (see, e.g. Lemma 9.4.A in \cite{Wol03}). 

As
$$\supp\mu\subset \bigcap_i\bigcup_{1\leq H\leq q_i^{s/2}}\mathcal{N}_{q_i^{-1}}\left(\frac{\Z}{H}\right)\cap[-\frac{1}{2}, \frac{1}{2}]$$
whose Hausdorff dimension is $\leq s$ by trivial coverings, the Fourier decay implies
$$\dH(\supp\mu)=s$$
and $\mu$ is a Salem measure.

We finish this section with the proof of Lemma \ref{stability-lemma}.
\begin{proof}[Proof of Lemma \ref{stability-lemma}]
Recall \eqref{Fourier-coefficient-F-i-mu}. By the prime number theorem, a trivial divisor bound, and the fact $\phi\in C^2$,
\begin{equation}\label{estimate-Fourier-coefficient-Fi-neq-0}|\widehat{F_i^\mu}(l)|\lesssim q_i^{-s/2}\log q_i\cdot(1+\log|l|)\cdot (1+\frac{|l|}{q_i})^{-2}, \ \forall\,l\in\Z\backslash\{0\}.\end{equation}
	As $\widehat{F_i^\mu}(0)=1$,
        \begin{equation}\label{sum-without-zero}\widehat{\psi F_i^\mu}(k)-\widehat{\psi}(k)=\sum_{l\in\mathbb{Z}}\widehat{\psi}(k-l)\widehat{F_i^\mu}(l)-\widehat{\psi}(k)
        =\sum_{l\neq0}\widehat{\psi}(k-l)\widehat{F_i^\mu}(l),\end{equation}
    and by \eqref{estimate-Fourier-coefficient-Fi-neq-0} its absolute value is 
        \begin{equation}\label{l-neq-0}\leq C q_i^{-s/2}\log q_i\cdot\sum_{l\neq 0}|\hat{\psi}(k-l)|\cdot(1+\log|l|)\cdot (1+\frac{|l|}{q_i})^{-2}.\end{equation}
        When $|k|\leq q_i$ and $|l|\leq 2q_i$, \eqref{l-neq-0} is
        $$\leq C q_i^{-s/2}(\log q_i)^2\cdot\|\hat{\psi}\|_{l^1}\leq C(\|\psi\|_{L^\infty}+\|\psi''\|_{L^\infty})\cdot q_i^{-s/2}(\log q_i)^2.$$
        When $|k|\leq q_i$ and $|l|\geq 2q_i$, we have $|k-l|\approx |l|$ and \eqref{l-neq-0} is
        $$\leq C q_i^{-s/2}\log q_i\cdot\sum_{|l|\geq q_i}|\hat{\psi}(l)|\log|l|\leq Cq_i^{-s/2}\log q_i\cdot\|\psi''\|_{L^\infty}\sum_{|l|\neq 0}|l|^{-2}\log|l|.$$
    From now we assume $|k|\geq q_i$ and write the sum \eqref{sum-without-zero} as
    $$\sum_{l\neq 0: |k-l|>|k|/2}+\sum_{l\neq 0: |k-l|\leq |k|/2}:=I+II.$$
    It is easy to estimate $I$: in this case $1\leq 2 |k|^{-1}|k-l|$, so, by $|\widehat{F_i^\mu}|\leq 1$,
    $$\sum_{l\neq0: |k-l|>|k|/2}|\widehat{\psi}(k-l)||\widehat{F_i^\mu}(l)|\leq \sum_{|k-l|>|k|/2}|\widehat{\psi}(k-l)|\leq  C|k|^{-s/2}\sum_{k-l\neq 0}|\hat{\psi}(k-l)||k-l|^{s/2}$$
    $$\leq C\|\psi''\|_{L^\infty}\cdot |k|^{-s/2}.$$
    For $II$, in this case $q_i/2\leq|k|/2\leq |l|\leq 3|k|/2$, so by \eqref{estimate-Fourier-coefficient-Fi-neq-0},
    $$\begin{aligned}
        &\sum_{l\neq0: |k-l|>|k|/2}|\widehat{\psi}(k-l)||\widehat{F_i^\mu}(l)|\\\leq &C\sum_{|k|/2\leq |l|\leq 3|k|/2}|\widehat{\psi}(k-l)| \cdot q_i^{-s/2}(\log|k|)^2\cdot(\frac{|l|}{q_i})^{-2}\\\leq &C \|\hat{\psi}\|_{l^1(\Z)}\cdot |k|^{-s/2}(\log|k|)^2\\\leq &C(\|\psi\|_{L^\infty}+\|\psi''\|_{L^\infty})\cdot |k|^{-s/2}(\log|k|)^2.
    \end{aligned}$$
\end{proof}

\section{Diophantine approximations: an auxiliary measure $\nu$}\label{sec-nu}
Even under such a careful construction, the spherical argument still does not apply to our $\mu$ directly. To prove Theorem \ref{thm-Diophantine} we need another measure $\nu$ to help, despite not knowing whether this auxiliary measure admits a Fourier frame.

To construct $\nu$, denote
$$\mathcal{P}^\nu_i:=\{q_i^{s/2}/h(i)\leq p\leq q_i^{s/2}/2, prime\},$$ where $h(i)$ is an auxiliary function. We shall see below that $h(i)=C_s\log q_i$, where $C_s>0$ only depends on $s$, is sufficient.

Then let $\Phi_{i,p}^\nu$ be a modified $p^{-1}$-periodization of $p^{-1}q_i\phi(q_ix)$,
\begin{equation}\label{def-Phi-i-p-nu}\begin{aligned}\Phi_{i,p}^\nu(x):=&\sum_{v\in\mathbb{Z}\backslash p\Z}p^{-1}q_i\phi(q_i(x-\frac{v}{p}))\\=&\sum_{n\in\Z}\hat{\phi}(pq_i^{-1}n) e^{2\pi i pn x}-p^{-1}\sum_{m\in\Z}\hat{\phi}(q_i^{-1}m)e^{2\pi i mx},\end{aligned}\end{equation}
and take
\begin{equation}\label{def-F-i-nu}F^\nu_i:=\frac{1}{\#\mathcal{P}_i^\nu}\sum_{p\in\mathcal{P}^\nu_i}\frac{p}{p-1}\Phi_{i,p}^\nu.\end{equation}
Notice $\widehat{F_i^\nu}(0)=1$ and for all integer $k\neq 0$,
\begin{equation}\label{Fourier-coefficient-F-i-nu}\widehat{F_i^\nu}(k)=
        (\#\mathcal{P}_i^\nu)^{-1}\left(\#\{p\in\mathcal{P}_i^\nu: k\in p\Z\}-\sum_{p\in\mathcal{P}_i^\nu: k\notin p\Z}\frac{1}{p-1}\right)\hat{\phi}(q_i^{-1}k).\end{equation}
Finally we take $\nu$ as the weak* limit of $\prod_{i=0}^n F^\nu_i$, with 
$$F^\nu_0(x):=2\phi(2x)=F^\mu_0(x)$$
as convention, which is a nonzero finite Borel measure supported on
$$\bigcap_i\bigcup_{p\in\mathcal{P}_i^\nu}\mathcal{N}_{q_i^{-1}}\left(\frac{\Z\backslash p\Z}{p}\right)\cap[-\frac{1}{2}, \frac{1}{2}].$$
As a trivial divisor bound shows
\begin{equation}\label{trivial-divisor-bound}|\#\{p\in\mathcal{P}_i^\nu:  k\in p\Z\}-\sum_{p\in\mathcal{P}_i^\nu:k\notin p\Z}\frac{1}{p-1}|\lesssim \frac{\log|k|}{\log (q_i^{s/2}/h(i))}+\frac{\#\mathcal{P}_i^\nu}{q_i^{s/2}/h(i)},\end{equation}
with $h(i)=C_s\log q_i$ in mind $|F^\nu_i|$ satisfies a better upper bound than $|F^\mu_i|$ in \eqref{estimate-Fourier-coefficient-Fi-neq-0}. So the stability lemma, Lemma \ref{stability-lemma}, also works on $\nu$ and then the weak* convergence follows in the same way as $\mu$. We omit details. It is a Salem measure as well.

Compared with $\mu$ from the previous section, there are two differences in the construction of $\nu$: the set of primes $\mathcal{P}_i^\nu$ and the exclusion of $p\Z$ in the definition of $\Phi_{i,p}^\nu$ in \eqref{def-Phi-i-p-nu}. In this section we explain the advantage of excluding $p\Z$. The relation between $\mu$ and $\nu$ will be carefully discussed in the next section. 

The exclusion of $p\Z$ ensures the following Frostman condition on $\nu$.
Although it makes $\hat{\nu}$ not necessarily nonnegative, we only need the positivity of $\hat{\mu}$.
\begin{lem}\label{Frostman-nu} For all $x\in\R$ and all $n\geq 1$,
	$$\nu(B(x,q_{n}^{-1}))\leq C^nq_{n}^{-s}\log q_{n}\prod_{i=1}^{n-1}q_i^{1-s}\log q_i\cdot \prod_{i=1}^n h(i),$$
	where the constant $C>0$ is independent in $x$ and $n$.
\end{lem}
In fact a more general Frostman condition holds on $\nu(B(x,r))$ for all $r>0$. One can see, for example, Section 4.3 in \cite{LL24+} for a proof on a similar construction. In this paper we only need the estimate on $B(x,q_n^{-1})$ which is more straightforward. So we give the proof for completeness.
\begin{proof}
	[Proof of Lemma \ref{Frostman-nu}]
	Recall the function $\phi_0$ we defined at the beginning of Section \ref{sec-mu}. Then it suffices to consider
	$$\int\sum_{v\in\Z}\phi_0(q_n(x+v-y))\,d\nu(y)=\sum_{k\in\Z}e^{2\pi i x\cdot k}q_n^{-1}\hat{\phi_0}(q_n^{-1}k)\hat{\nu}(k).$$
	We have pointed out that the stability lemma, Lemma \ref{stability-lemma}, also works on $\nu$. As $\nu$ is the weak* limit of $\prod F_i^\nu$ and $q_i$ is a rapidly increasing sequence, we have
	\begin{equation}\label{reduction-to-finite-product}|\hat{\nu}(k)-\widehat{\prod_{0\leq i\leq n}F_i^\nu}(k)|\leq \sum_{m=n+1}^\infty |\widehat{\prod_{0\leq i\leq m}F_i^\nu}(k)-\widehat{\prod_{0\leq i\leq m-1}F_i^\nu}(k)|\lesssim_\epsilon
		q_{n+1}^{-s/2+\epsilon},\ \forall\,k\in\Z,\end{equation}
	where the implicit constant is independent in $k$ and $n$. Of course this estimate also holds on $\mu$.
	
	If we choose, for example $\epsilon=s/4$, it follows that
	$$\left|\int\sum_{v\in\Z}\phi_0(q_n(x+v-y))\,d\nu(y)-\int\sum_{v\in\Z}\phi_0(q_n(x+v-y))\prod_{0\leq i\leq n}F_i^\nu(y)\,dy\right|$$ $$\lesssim q_{n+1}^{-s/4}\sum_{k\in\Z}q_n^{-1}\hat{\phi_0}(q_n^{-1}k)\lesssim q_{n+1}^{-s/4},$$
	much better than the expected upper bound due to the rapid increase of $q_i$. Therefore it suffices to consider
	$$\int\sum_{v\in\Z}\phi_0(q_n(x+v-y))\prod_{0\leq i\leq n}F_i^\nu(y)\,dy\lesssim q_n^{-1}\prod_{0\leq i\leq n}\|F_i^\nu\|_{L^\infty}.$$
	
	Here comes the advantage of excluding $p\Z$. Recall in the construction of $\mu$, we point out that for each fixed $p\in \mathcal{P}_i^\mu$, translated $\phi$ have disjoint support in definition of of $\Phi_{i,p}^\mu$. Thanks to the exclusion of $p\Z$, all translated $\phi$ have disjoint support in the definition  of $F_i^\nu$, even for those from different $p$. This is because for all primes $p, p'\in \mathcal{P}_i^\nu$ and integers $m\notin p\Z$, $m'\notin p'\Z$,
    \begin{equation}\label{separation}\left|\frac{m}{p}-\frac{m'}{p'}\right|=\frac{|mp'-m'p|}{pp'}> 4q_i^{-s}\geq 4q_i^{-1}, \text{ given }(m,p)\neq(m',p').\end{equation}
	
	As a consequence,
	$$\|F_i^\nu\|_{L^\infty}\leq  (\#\mathcal{P}_i^\nu)^{-1}\cdot\sup_{p\in\mathcal{P}_i^\nu}\frac{q_i}{p-1}\|\phi\|_{L^\infty}\lesssim q_i^{1-s}\cdot \log q_i\cdot h(i)$$
	which is a better estimate than that of $F_i^\mu$ in \eqref{F-i-mu-L-infty}.
	
	Overall,
	$$\nu(B(x,q_n^{-1})\lesssim q_n^{-1}\prod_{0\leq i\leq n}\|F_i^\nu\|_{L^\infty}\leq C^nq_n^{-s}\log q_n\prod_{i=1}^{n-1}q_i^{1-s}\log q_i\cdot \prod_{i=1}^n h(i),$$
	as desired.
\end{proof}

\section{Diophantine approximations: relation between $\mu$ and $\nu$}\label{sec-mu-nu}
In this section we make careful comparison between the target measure $\mu$ and the auxiliary measure $\nu$. There are two key lemmas.

\begin{lem}
\label{lem-absolute-continuity}
$\nu\ll\mu$ with $\frac{d\nu}{d\mu}\in L^\infty(\mu)$.	
\end{lem}

\begin{lem}
	\label{lem-averaging-frequency}Suppose $|k|>2q_n, |l|< q_n/2$ are integers. Then
	$$|\hat{\nu}(k+l)|\leq C\hat{\mu}(k)+C_\epsilon (1+|k|)^{-1+\epsilon},$$ where $C, C_\epsilon>0$ are independent in $n$, $k$ and $l$.
\end{lem}

\begin{proof}
	[Proof of Lemma \ref{lem-absolute-continuity}]
	It suffices to prove $\nu(E)\leq C\mu(E)$ for all Borel sets $E$. As
	$$\mu=\lim_{n\rightarrow\infty}\prod_{i=0}^n F_i^\mu\quad \text{and}\quad \nu=\lim_{n\rightarrow\infty}\prod_{i=0}^n F_i^\nu,$$
	we compare $F_i^\mu$ defined in \eqref{def-F-i-mu} and $F_i^\nu$ defined in \eqref{def-F-i-nu}. From their construction,
	$$\begin{aligned}F_i^\nu(x)= &\frac{1}{\#\mathcal{P}_i^\nu}\sum_{p\in\mathcal{P}^\nu_i}\frac{p}{p-1}\sum_{v\in\mathbb{Z}\backslash p\Z}p^{-1}q_i\phi(q_i(x-\frac{v}{p}))\\\leq &\frac{q_i^{s/2}/h(i)}{q_i^{s/2}/h(i)-1}\cdot \frac{\#\mathcal{P}_i^\mu}{\#\mathcal{P}_i^\nu}\cdot\frac{1}{\#\mathcal{P}_i^\mu}\sum_{p\in\mathcal{P}^\nu_i}\sum_{v\in\mathbb{Z}\backslash p\Z}p^{-1}q_i\phi(q_i(x-\frac{v}{p}))\\\leq &\frac{q_i^{s/2}/h(i)}{q_i^{s/2}/h(i)-1}\cdot \frac{\#\mathcal{P}_i^\mu}{\#\mathcal{P}_i^\nu}\cdot\frac{1}{\#\mathcal{P}_i^\mu}\sum_{p\in\mathcal{P}^\mu_i}\sum_{v\in\mathbb{Z}}p^{-1}q_i\phi(q_i(x-\frac{v}{p}))\\=&\frac{q_i^{s/2}/h(i)}{q_i^{s/2}/h(i)-1}\cdot \frac{\#\mathcal{P}_i^\mu}{\#\mathcal{P}_i^\nu}\cdot F_i^\mu(x).
	\end{aligned}$$
	
	If $h(i)$ is chosen properly, say $h(i)=C_s\log q_i$, then both infinite products
	$$\prod \frac{q_i^{s/2}/h(i)-1}{q_i^{s/2}/h(i)}=\prod (1-q_i^{-s/2}h(i))$$
	and
	$$\prod \frac{\#\mathcal{P}_i^\nu}{\#\mathcal{P}_i^\mu}=\prod\left(1- \frac{1+\#\{p\leq q_i^{s/2}/h(i)\}}{1+\#\{p\leq q_i^{s/2}/2\}}\right)$$
	are convergent. Hence
	$$\lim_{n\rightarrow\infty}\prod_{i=0}^n F_i^\mu\leq C\lim_{n\rightarrow\infty}\prod_{i=0}^n F_i^\nu$$
	in the weak* sense, as desired.
\end{proof}

\begin{proof}
	[Proof of Lemma \ref{lem-averaging-frequency}]
	Suppose $2q_{n_0-1}< |k|\leq 2q_{n_0}$, $n_0\geq n+1$. By the estimate \eqref{reduction-to-finite-product}, up to a negligible error $\hat{\nu}(k+l)$ equals
$$\widehat{\prod_{0\leq i\leq n_0} F_i^\nu}(k+l)=\sum_{k_0+\cdots+k_{n_0}=k+l} \widehat{F_0^\nu}(k_0)\cdots\widehat{F_{n_0}^\nu}(k_{n_0})$$ 
$$=\sum_{\substack{m_0+\cdots+m_{n_0}=k\\ m_{n_0}+l\neq 0}} \left(\prod_{0=1}^{n_0-1}\widehat{F_i^\nu}(m_i)\right)\widehat{F_{n_0}^\nu}(m_{n_0}+l)+\sum_{m_0+\cdots+m_{n_0-1}=k} \left(\prod_{i=0}^{n_0-2}\widehat{F_i^\nu}(m_i)\right)\widehat{F_{n_0-1}^\nu}(m_{n_0-1}+l)$$
$$:=I+II.$$

Recall $\widehat{F_i^\mu}, \widehat{F_i^\nu}$ have been expressed in \eqref{Fourier-coefficient-F-i-mu}, \eqref{Fourier-coefficient-F-i-nu}. Also $\|\widehat{F_i^\mu}\|_{l^1(\Z)}=F_i^\mu(0)\leq \|F_i^\mu\|_{L^\infty}$ because $\widehat{F_i^\mu}\geq 0$ on $\Z$.

Fix $i\geq 1$. With $h(i)=C_s\log q_i$ and $C_s>0$ properly chosen, for all integer $m\neq 0$,
$$\begin{aligned}|\#\{p\in\mathcal{P}_i^\nu: m\in p\Z\}-\sum_{p\in\mathcal{P}_i^\nu: k\notin p\Z}\frac{1}{p-1}| \leq &\#\{p\in\mathcal{P}_i^\nu: m\in p\Z\}+ \frac{\#\mathcal{P}_i^\nu}{q_i^{s/2}/h(i)-1}\\\leq &\#\{p\in\mathcal{P}_i^\nu: m\in p\Z\}+1\\\leq &\#\{p\in\mathcal{P}_i^\mu: m\in p\Z\}.\end{aligned}$$
Therefore, for each $i\geq 1$ and $m\in\Z\backslash\{0\}$,
\begin{equation}\label{first-comparison}|\widehat{F_i^\nu}(m)|\leq \frac{\#\mathcal{P}_i^\mu}{\#\mathcal{P}_i^\nu} \cdot (\#\mathcal{P}_i^\mu)^{-1}\cdot\#\{p\in\mathcal{P}_i^\mu: m\in p\Z\}\cdot\hat{\phi}(q_i^{-1}m)=\frac{\#\mathcal{P}_i^\mu}{\#\mathcal{P}_i^\nu} \cdot\widehat{F_i^\mu}(m).\end{equation}
When $m=0$ this estimate trivially holds.

We estimate $I$ first.

For $|m_{n_0}|>q_{n_0}^{10}$, by the prime number theorem, \eqref{trivial-divisor-bound}, and the Fourier decay $\hat{\phi}(\xi)\lesssim (1+|\xi|)^{-2}$ as $\phi\in C^2(-1,1)$,
$$|\widehat{F_n^\nu}(m_{n_0}+l)|\lesssim q_{n_0}^{-s/2}\log q_{n_0} \log |m_{n_0}| \cdot \left(\frac{|m_{n_0}|}{q_{n_0}}\right)^{-2}\lesssim q_{n_0}^{-10},$$
where the implicit constant is independent in $n_0$, $m_{n_0}$, $l$, and therefore
$$I\lesssim q_{n_0}^{-10}\prod_{i=0}^{n_0-1}\|\widehat{F_i^\nu}\|_{l^1(\Z)}\leq q_{n_0}^{-10}\prod_{i\geq 1}\frac{\#\mathcal{P}_i^\mu}{\#\mathcal{P}_i^\nu}\cdot \prod_{i=0}^{n_0-1}\|\widehat{F_i^\mu}\|_{l^1(\Z)}\lesssim q_{n_0}^{-10}\prod_{i=0}^{n_0-1}\|F_i^\mu\|_{L^\infty},$$
which is $\lesssim q_{n_0}^{-9}$, a negligible error, by the convergence of $\prod\frac{\#\mathcal{P}_i^\mu}{\#\mathcal{P}_i^\nu}$, the estimate \eqref{F-i-mu-L-infty} on $\|F_i^\mu\|_{L^\infty}$, and the rapid increase of $q_i$.

For $|m_{n_0}|<q_{n_0}^{10}$ and $m_{n_0}\neq 0$, the trivial divisor bound \eqref{trivial-divisor-bound} on $m_{n_0}+l$ is $\lesssim_s 1$, given $h(i)=C_s\log q_i$. Therefore
$$|\widehat{F_{n_0}^\nu}(m_{n_0}+l)|\lesssim (\#\mathcal{P}_{n_0}^\nu)^{-1}\cdot \hat{\phi}(q_{n_0}^{-1}(m_{n_0}+l))\lesssim (\#\mathcal{P}_{n_0}^\nu)^{-1}\cdot \hat{\phi}(q_{n_0}^{-1}m_{n_0}),$$
where the last inequality follows from the last property of $\phi$ listed in Section \ref{sec-mu} and the fact $|q_{n_0}^{-1}l|\leq 1$. As $\#\mathcal{P}_{n_0}^\nu\approx \#\mathcal{P}_{n_0}^\mu$ and $\#\{p\in\mathcal{P}_{n_0}^\mu: m\in p\Z\}\geq 1$ for all $m\in\Z$, it follows that
\begin{equation}\label{second-comparison}|\widehat{F_{n_0}^\nu}(m_{n_0}+l)|\lesssim (\#\mathcal{P}_{n_0}^\mu)^{-1} \cdot \#\{p\in\mathcal{P}_{n_0}^\mu: m_{n_0}\in p\Z\}\cdot\hat{\phi}(q_i^{-1}m_{n_0})=\widehat{F_n^\mu}(m_{n_0}),\end{equation}
where the implicit constant is independent in $n_0$, $m_{n_0}$, $l$. When $m_{n_0}=0$ this estimate trivially holds as $|\widehat{F_{n_0}^\nu}|\leq 1$.

Putting \eqref{first-comparison}\eqref{second-comparison} together, up to an error $q_{n_0}^{-10}\lesssim |k|^{-1}$,
$$I\lesssim \prod_{i\geq 1}\frac{\#\mathcal{P}_i^\mu}{\#\mathcal{P}_i^\nu}\cdot\sum_{m_0+\cdots+m_{n_0}=k} \prod_{i=0}^{n_0}\widehat{F_i^\mu}(m_i)\lesssim \sum_{m_0+\cdots+m_{n_0}=k} \prod_{i=0}^{n_0}\widehat{F_i^\mu}(m_i).$$
As \eqref{reduction-to-finite-product} also holds on $\mu$, the right hand side equals $\hat{\mu}(k)$ up to a negligible error.

Now we estimate $II$. 

The first case is $|m_{n_0-1}+l|>\max\{q_{n_0-1}^{10}, |k|/2\}$. By the prime number theorem, \eqref{trivial-divisor-bound}, and the decay of $\hat{\phi}$,
$$|\widehat{F_{n_0-1}^\nu}(m_{n_0-1}+l)|\lesssim q_{n_0}^{-s/2}\log q_{n_0} \log |m_{n_0}+l| \cdot \left(\frac{|m_{n_0}+l|}{q_{n_0}}\right)^{-2}\lesssim |k|^{-1}q_{n_0-1}^{-5},$$
where the implicit constant is independent in $n_0$, $m_{n_0-1}$, $l$, and therefore
$$II\lesssim |k|^{-1}q_{n_0-1}^{-5}\prod_{i=0}^{n_0-2}\|\widehat{F_i^\nu}\|_{l^1(\Z)}\leq |k|^{-1}q_{n_0-1}^{-5}\prod_{i\geq 1}\frac{\#\mathcal{P}_i^\mu}{\#\mathcal{P}_i^\nu}\cdot \prod_{i=0}^{n_0-2}\|\widehat{F_i^\mu}\|_{l^1(\Z)}\lesssim |k|^{-1}q_{n_0-1}^{-5}\prod_{i=0}^{n_0-2}\|F_i^\mu\|_{L^\infty},$$
which is $\lesssim |k|^{-1}$ by the convergence of $\prod\frac{\#\mathcal{P}_i^\mu}{\#\mathcal{P}_i^\nu}$, the estimate \eqref{F-i-mu-L-infty} on $\|F_i^\mu\|_{L^\infty}$, and the rapid increase of $q_i$.

The second case is $|m_{n_0-1}+l|\leq |k|/2$. In this case $|m_{n_0-1}|\leq 3|k|/4$, which means there must exist $1\leq i_0\leq n_0-2$ such that 
$$|m_{i_0}|>(4n_0)^{-1}|k|>(2n_0)^{-1}q_{n_0-1}.$$ 
So by the prime number theorem, \eqref{trivial-divisor-bound}, and the decay of $\hat{\phi}$,
$$|\widehat{F_{i_0}^\nu}(m_{i_0})|\lesssim q_{i_0}^{-s/2}\log q_{i_0} \log |m_{i_0}| \cdot \left(\frac{|m_{i_0}|}{q_{i_0}}\right)^{-2}\lesssim_\epsilon |k|^{-1+2\epsilon}q_{n_0-1}^{-1-\epsilon},$$
and therefore
$$II\lesssim_\epsilon |k|^{-1+2\epsilon}q_{n_0-1}^{-1-\epsilon}\prod_{i=0}^{n_0-1}\|\widehat{F_i^\nu}\|_{l^1(\Z)}\leq |k|^{-1+2\epsilon}q_{n_0-1}^{-1-\epsilon}\prod_{i\geq 1}\frac{\#\mathcal{P}_i^\mu}{\#\mathcal{P}_i^\nu}\cdot \prod_{i=0}^{n_0-1}\|\widehat{F_i^\mu}\|_{l^1(\Z)}$$ $$\lesssim |k|^{-1+2\epsilon}q_{n_0-1}^{-1-\epsilon}\prod_{i=0}^{n_0-1}\|F_i^\mu\|_{L^\infty},$$
which is $\lesssim_\epsilon |k|^{-1+2\epsilon}$ by the convergence of $\prod\frac{\#\mathcal{P}_i^\mu}{\#\mathcal{P}_i^\nu}$, the estimate \eqref{F-i-mu-L-infty} on $\|F_i^\mu\|_{L^\infty}$, and the rapid increase of $q_i$.

The last case is $|m_{n_0-1}+l|\leq q_{n_0-1}^{10}$. In this case the trivial divisor bound \eqref{trivial-divisor-bound} on $m_{n_0-1}+l$ is $\lesssim_s 1$, given $h(i)=C_s\log q_i$. Therefore for $m_{n_0-1}\neq 0$,
$$|\widehat{F_{n_0-1}^\nu}(m_{n_0-1}+l)|\lesssim (\#\mathcal{P}_{n_0-1}^\nu)^{-1}\cdot \hat{\phi}(q_{n_0-1}^{-1}(m_{n_0-1}+l))\lesssim (\#\mathcal{P}_{n_0-1}^\nu)^{-1}\cdot \hat{\phi}(q_{n_0-1}^{-1}m_{n_0-1}),$$
where the last inequality follows from the last property of $\phi$ listed in Section \ref{sec-mu} and the fact $|q_{n_0-1}^{-1}l|\leq 1$. As $\#\mathcal{P}_{n_0-1}^\nu\approx \#\mathcal{P}_{n_0-1}^\mu$ and $\#\{p\in\mathcal{P}_{n_0-1}^\mu: m\in p\Z\}\geq 1$ for all $m\in\Z$, it follows that
\begin{equation}\label{third-comparison}|\widehat{F_{n_0-1}^\nu}(m_{n_0-1}+l)|\lesssim (\#\mathcal{P}_{n_0-1}^\mu)^{-1} \cdot \#\{p\in\mathcal{P}_{n_0-1}^\mu: p|m_{n_0-1}\}\cdot\hat{\phi}(q_{n_0-1}^{-1}m_{n_0-1})=\widehat{F_{n_0-1}^\mu}(m_{n_0-1}),\end{equation}
where the implicit constant is independent in $n_0$, $m_{n_0-1}$, $l$. When $m_{n_0-1}=0$ this estimate trivially holds as $|\widehat{F_{n_0-1}^\nu}|\leq 1$. 

Putting \eqref{first-comparison}\eqref{third-comparison} together, in this case 
$$II\lesssim \prod_{i\geq 1}\frac{\#\mathcal{P}_i^\mu}{\#\mathcal{P}_i^\nu}\cdot\sum_{m_0+\cdots+m_{n_0-1}=k} \prod_{i=0}^{n_0-1}\widehat{F_i^\mu}(m_i)\lesssim\sum_{m_0+\cdots+m_{n_0}=k} \prod_{i=0}^{n_0}\widehat{F_i^\mu}(m_i).$$
As \eqref{reduction-to-finite-product} also holds on $\mu$, the right hand side equals $\hat{\mu}(k)$ up to a negligible error.
\end{proof}

\section{Diophantine approximations: no Fourier frames on $\mu$}\label{sec-final-proof}
In this section we finish the proof of Theorem \ref{thm-Diophantine} with $\mu$ constructed in Section \ref{sec-mu}.

Suppose $\Lambda\subset\R$ is a Fourier frame for $\mu$, namely
$$A\|f\|_{L^2(\mu)}\leq\sum_{\lambda\in\Lambda}|\widehat{f\,d\mu}(\lambda)|^2\leq B\|f\|_{L^2(\mu)},\ \forall\,f\in L^2(\mu).$$
By considering $f(x)e^{-2\pi ix\xi}$, it becomes
\begin{equation}\label{frame-mu}A\|f\|_{L^2(\mu)}\leq\sum_{\lambda\in\Lambda}|\hat{f\,d\mu}(\lambda+\xi)|^2\leq B\|f\|_{L^2(\mu)},\ \forall\,\xi\in\R.\end{equation}
In particular,
$$	\sum_{\lambda\in\Lambda}|\hat{\mu}(\lambda+\xi)|^2\approx 1$$
uniformly in $\xi\in\R$. It implies a basic property of $\Lambda$: for all $r>0$,
\begin{equation}\label{frame-non-concentration}\sup_x\#\{\Lambda\cap B(x,r)\}<\infty.\end{equation}
To see this, there exists $\delta>0$ such that $|\hat{\mu}|\geq C>0$ on $[-\delta, \delta]$, then for all $x\in\R$,
$$B\geq \sum_{\lambda\in\Lambda}|\hat{\mu}(\lambda-x)|^2\geq C\#\{\Lambda\cap B(x,\delta)\}.$$
As a consequence of \eqref{frame-non-concentration},
\begin{equation}\label{convergence-lambda}\sum_{\lambda\in\Lambda\backslash\{0\}} |\lambda|^{-\alpha}<\infty,\ \forall\alpha>1.\end{equation}

By Lemma \ref{lem-absolute-continuity}, $\frac{d\nu}{d\mu}$ is a compactly supported $L^\infty(\mu)$ function, so $\frac{d\nu}{d\mu}\in L^2(\mu)$, and $\|\frac{d\nu}{d\mu}\|_{L^2(\mu)}\neq 0$ because $\nu$ is a nontrivial measure. Also, if we take 
$$\psi(x):=\phi_0(2x)\in C_0^\infty(-\frac{1}{2},\frac{1}{2})$$
with $\phi_0$ defined at the beginning of Section \ref{sec-mu}, then $\|\psi\|_{L^2(\mu)}, \|\psi\|_{L^2(\nu)}\neq 0$. This is because by \eqref{reduction-to-finite-product}, up to a negligible error, both
$$\mu(B(0,q_n^{-1}))=\int_{-q_n^{-1}}^{q_n^{-1}}\prod_{i=0}^nF_i^\mu>0 \text{ and } \nu(B(0,q_n^{-1}))=\int_{-q_n^{-1}}^{q_n^{-1}}\prod_{i=0}^nF_i^\nu>0.$$
We omit details. 

Consequently, by plugging
$$f(x)=\psi(x)\frac{d\nu}{d\mu}(x)e^{-2\pi ix\xi}$$
into \eqref{frame-mu}, we obtain
$$A'\leq \sum_{\lambda\in\Lambda}|\widehat{\psi\,d\nu}(\lambda+\xi)|^2\leq B'$$
for some constants $0<A'\leq B'<\infty$. By taking an average over $\xi\in [-\frac{q_n}{2}, \frac{q_n}{2}]$, it becomes
\begin{equation}
	\label{sum-avg-psi-nu}
	A'\leq \sum_{\lambda\in\Lambda}q_n^{-1}\int_{-q_n/2}^{q_n/2}|\widehat{\psi\,d\nu}(\lambda+\xi)|^2\,d\xi\leq B',\ \forall\,n.
\end{equation}	
As
\begin{equation}\label{hat-lambda+xi}
	\widehat{\psi\,d\nu}(\lambda+\xi)=\int_{-\frac{1}{2}}^{\frac{1}{2}} e^{-2\pi i (\lambda+\xi) x}\phi_0(2x)\,d\nu(x)=\frac{1}{2}\sum_{k\in\Z}\hat{\phi_0}(\frac{\lambda+\xi-k}{2})\hat{\nu}(k),
\end{equation}
by splitting $[-\frac{q_n}{2}, \frac{q_n}{2}]$ into unit intervals it follows that
\begin{equation}\label{avg-psi-nu}
	\begin{aligned}q_n^{-1}\int_{-q_n/2}^{q_n/2}|\widehat{\psi\,d\nu}(\lambda+\xi)|^2\,d\xi\leq & q_n^{-1}\sum_{|l|\leq q_n/2}\int_{-1}^1\left(\sum_{k\in\Z}\hat{\phi_0}(\frac{\lambda+t+l-k}{2})|\hat{\nu}(k)|\right)^2\,dt\\=&q_n^{-1}\sum_{|l|\leq q_n/2}\int_{-1}^1\left(\sum_{k\in\Z}\hat{\phi_0}(\frac{\lambda+t-k}{2})|\hat{\nu}(k+l)|\right)^2\,dt.\end{aligned}
\end{equation}

If $|\lambda|>4q_n$ and $|k|>|\lambda|/2$, by Lemma \ref{lem-averaging-frequency} the last line of \eqref{avg-psi-nu} is
\begin{equation}\label{ave-psi-nu-main-term}\begin{aligned}
&\lesssim q_n^{-1}\sum_{|l|\leq q_n/2}\int_{-1}^1\left(\sum_{|k|>|\lambda|/2}\hat{\phi_0}(\frac{\lambda+t-k}{2})(\hat{\mu}(k)+|k|^{-0.9})\right)^2\,dt\\&\lesssim \int_{-1}^1\left(\frac{1}{2}\sum_{k\in\Z}\hat{\phi_0}(\frac{\lambda+t-k}{2})\hat{\mu}(k)\right)^2\,dt+|\lambda|^{-1.8}\int_{-1}^1\left(\sum_{k\in\Z}\hat{\phi_0}(\frac{\lambda+t-k}{2})\right)^2\,dt\\&=\int_{-1}^1|\widehat{\psi\,d\mu}(\lambda+t)|^2\,dt +O(|\lambda|^{-1.8}).
\end{aligned}\end{equation}

If $|\lambda|>4q_n$ and $|k|\leq |\lambda|/2$, we have $|\lambda+t-k|>|\lambda|/4$. Then due to $|\hat{\nu}|\leq 1$ and the fast decay of $\hat{\phi_0}$, the last line of \eqref{avg-psi-nu} is
\begin{equation}\label{ave-psi-nu-error-term}\lesssim_N|\lambda|^{-N}.\end{equation}
By \eqref{avg-psi-nu}\eqref{ave-psi-nu-main-term}\eqref{ave-psi-nu-error-term}, one can take the sum over $|\lambda|>4q_n$ to obtain
$$\sum_{|\lambda|>4q_n}q_n^{-1}\int_{-q_n/2}^{q_n/2}|\widehat{\psi\,d\nu}(\lambda+\xi)|^2\,d\xi\lesssim \sum_{|\lambda|>4q_n}\int_{-1}^1|\widehat{\psi\,d\mu}(\lambda+t)|^2\,dt +|\lambda|^{-1.8},$$
which goes to $0$ as $n\rightarrow\infty$ due to the convergence of 
$$\sum_{\lambda\in\Lambda\backslash\{0\}}\int_{-1}^1|\widehat{\psi\,d\mu}(\lambda+t)|^2\,dt +|\lambda|^{-1.8}$$
discussed in \eqref{frame-mu}\eqref{convergence-lambda} above. This means, if we split the sum in \eqref{sum-avg-psi-nu} into 
$$\sum_{|\lambda|\leq 4q_n}+\sum_{|\lambda|>4q_n},$$
then one can drop the second term when $n$ is large to conclude
\begin{equation}\label{reduction-to-counting}\begin{aligned}A'/2\leq &\sum_{|\lambda|\leq 4q_n}q_n^{-1}\int_{-q_n/2}^{q_n/2}|\widehat{\psi\,d\nu}(\lambda+\xi)|^2\,d\xi\\\leq &\#\{\Lambda\cap[-4q_n, 4q_n]\}\cdot q_n^{-1}\int_{-5q_n}^{5q_n}|\widehat{\psi\,d\nu}(\xi)|^2\,d\xi.\end{aligned}\end{equation}
Recall $\phi_0$ is defined at the beginning of Section \ref{sec-mu}. Then it is a standard argument that
$$q_n^{-1}\int_{-5q_n}^{5q_n}|\widehat{\psi\,d\nu}(\xi)|^2\,d\xi\lesssim q_n^{-1}\int|\widehat{\psi\,d\nu}(\xi)|^2\hat{\phi_0}(cq_n^{-1}\xi)\,d\xi$$
$$=C\iint\phi_0(c^{-1}q_n(x-y))\phi_0(2x)\phi_0(2y)d\nu(x)d\nu(y)\lesssim \|\phi_0\|_{L^\infty}^3\cdot \sup_x\nu(B(x,cq_n)^{-1})),$$
which is, by Lemma \ref{Frostman-nu},
$$\leq C^nq_{n}^{-s}\log q_{n}\prod_{i=1}^{n-1}q_i^{1-s}\log q_i\cdot \prod_{i=1}^n h(i).$$
Plugging this into \eqref{reduction-to-counting}, we obtain the counting estimate
\begin{equation}\label{counting-estimate}\#\{|\lambda|\leq 4q_n\}\geq \left(C^nq_n^{-s}\log q_n\prod_{i=1}^{n-1}q_i^{1-s}\log q_i\right)^{-1}\cdot\left(\prod_{i=1}^n h(i)\right)^{-1}.\end{equation}

On the other hand, as \eqref{hat-lambda+xi} also holds on $\mu$ and $\hat{\phi_0}\geq 0, \hat{\mu}(k)>0$, for all $\lambda\in\R$,
\begin{equation}\label{lower-bound-mu-1}\begin{aligned}
\int_{-1}^1|\widehat{\psi\,d\mu}(\lambda+t)|^2\,dt=&\int_{-1}^1\left|\frac{1}{2}\sum_{k\in\Z}\hat{\phi_0}(\frac{\lambda+t-k}{2})\hat{\mu}(k)\right|^2dt,
\\\geq &|\hat{\mu}([\lambda])|^2 \int_{-1}^1\left|\frac{1}{2}\hat{\phi_0}(\frac{\lambda+t-[\lambda]}{2})\right|^2dt\\\gtrsim&|\hat{\mu}([\lambda])|^2,\end{aligned}\end{equation}
where $[\cdot]$ denotes the integer part. 

Recall $\mu$ is the weak* limit of $\prod F_i^\mu$ and $\widehat{F_i^\mu}>0$ on $\Z$, so for every integer $|k|\leq 4q_n$,
$$\hat{\mu}(k)=\sum_{m=0}^\infty\sum_{\substack{k_0+\cdots+k_m=k\\ k_m\neq 0}} \prod_{i=0}^m\widehat{F_i^\mu}(k_i)\geq \sum_{k_1+\cdots+k_n=k} \prod_{i=1}^n\widehat{F_i^\nu}(k_i)\geq \sum_{\substack{k_1+\cdots+k_n=k\\|k_i|\leq q_i, 1\leq i\leq n-1}} \prod_{i=1}^n\widehat{F_i^\nu}(k_i).$$
As $|k_i|\leq q_i$ for all $i\leq n-1$,  the last term $|k_n|\leq 5q_n$. Therefore by the expression of $\widehat{F_i^\mu}$ in \eqref{Fourier-coefficient-F-i-mu}, we have
$$\widehat{F_i^\mu}(k_i)\geq c q_i^{-s/2}\log q_i,\ i=1,2,\dots,n,$$
and therefore
\begin{equation}\label{lower-bound-hat-mu-2}\hat{\mu}(k)\geq c q_n^{-s/2}\log q_n\cdot\prod_{i=1}^{n-1}\sum_{|k_i|\leq q_i}\widehat{F_i^\mu}(k_i)\geq c^n q_{n-1}^{-s/2}\log q_{n-1}\cdot\prod_{i=1}^{n-2}q_i^{1-s/2}\log q_i.\end{equation}

Now we are ready for the contradiction: by plugging \eqref{lower-bound-mu-1} into the upper bound of \eqref{frame-mu} with $f=\psi$, we have
$$2B\|\psi\|_{L^2(\mu)}^2\geq \sum_{\lambda\in\Lambda}\int_{-1}^1|\widehat{\psi\,d\mu}(\lambda+t)|^2\,dt\gtrsim\sum_{|\lambda|\leq 4q_n} |\hat{\mu}([\lambda])|^2,$$
which is, by \eqref{lower-bound-hat-mu-2} and \eqref{counting-estimate},
$$\geq \left(c^nq_n^{-s/2}\log q_n\cdot\prod_{i=1}^{n-1}q_i^{1-s/2}\log q_i\right)^2\cdot \left(C^nq_n^{-s}\log q_n\prod_{i=1}^{n-1}q_i^{1-s}\log q_i\right)^{-1}\left(\prod_{i=1}^n h(i)\right)^{-1}$$ $$= c^n \log q_n\prod_{i=1}^{n-1} q_i\log q_i/\prod_{i=1}^n h(i)\rightarrow\infty,$$
if $h(i)$ is chosen properly, say $h(i)=C_s\log q_i$.

\section{Spherical arcs are Salem}\label{sec-arc}
In this section we explain why the weighted arc mentioned at the end of introduction is Salem, that is
$$\int_{-\frac{1}{2}}^{\frac{1}{2}} e^{-2\pi i (x, \sqrt{1-x^2})\cdot(\xi_1, \xi_2)}\,dx\lesssim |\xi|^{-1/2}.$$
Although this is not a standard stationary phase estimate which requires the integrand to be $C_0^\infty$, the stationary phase technique still works (see e.g. \cite{Sogge17}, Chapter 1). For readers' convenience we sketch the proof for this special case.

Write $\xi=|\xi|(\cos\theta_0, \sin\theta_0)$. After changing variables
$$\int_{-\frac{1}{2}}^{\frac{1}{2}} e^{-2\pi i (x, \sqrt{1-x^2})\cdot(\xi_1, \xi_2)}\,dx=\int_{\theta_1}^{\theta_2} e^{-2\pi i (\cos\theta, \sin\theta)\cdot(\xi_1, \xi_2)}\,\sin\theta d\theta=\int_{\theta_1}^{\theta_2} e^{-2\pi i|\xi| \cos(\theta-\theta_0)}\,\sin\theta d\theta.$$
Decompose the interval $[\theta_1, \theta_2]$ into
$$I:=\{\theta\in [\theta_1, \theta_2]: \dist(\theta_0-\theta, 2\pi\Z)\geq |\xi|^{-1/2}\} $$
and its complement (if nonempty). Trivially
$$\left|\int_{[\theta_1, \theta_2]\backslash I} e^{-2\pi i|\xi| \cos(\theta-\theta_0)}\,\sin\theta d\theta\right|\lesssim |\xi|^{-1/2}.$$
 Then, as $I$ consists one or two intervals, by integration by parts twice
$$\begin{aligned}\left|\int_{I} e^{-2\pi i|\xi| \cos(\theta-\theta_0)}\,\sin\theta d\theta\right|\lesssim &|\xi|^{-1}\int_I \frac{1}{\sin^2(\theta-\theta_0)} d\theta+ |\xi|^{-1/2}\\\lesssim &|\xi|^{-1}\int_{|\xi|^{-1/2}}^\infty t^{-2} dt+ |\xi|^{-1/2}\\\lesssim & |\xi|^{-1/2},\end{aligned}$$
as desired, where $|\xi|^{-1/2}$ in the first line comes from the endpoint terms in the integration by parts.

Here the estimate is as good as the standard stationary phase because the endpoint terms coming from integration by parts do not dominate the integral term, which may not be the case in higher dimensions.

\bibliographystyle{abbrv}
\bibliography{mybibtex.bib}

\end{document}